\documentclass{amsart}
\usepackage{graphicx}
\usepackage{fancyvrb}
\usepackage{amssymb}
\usepackage{longtable}
\usepackage[pdfauthor={Richard J. Mathar},pdfkeywords={small graphs, multigraphs, graph enumeration},pdfsubject={Mathematics, Graph Theory},colorlinks=true,citecolor=blue]{hyperref}

\usepackage{amsmath}

\newtheorem{exa}{Example}
\newtheorem{defn}{Definition}
\newtheorem{remark}{Remark}

\begin{document}

\title[Small Unlabeled Graphs]{Statistics on Small Graphs}

\author{Richard J. Mathar}
\urladdr{http://www.mpia.de/~mathar}
\address{Max-Planck Institute of Astronomy, K\"onigstuhl 17, 69117 Heidelberg, Germany}


\date{\today}

\begin{abstract}
We create the unlabeled or vertex-labeled
graphs with up to 10 edges and up to 10 vertices
and classify them by a set of standard properties: directed
or not, vertex-labeled or not, connectivity, presence of isolated vertices, presence
of multiedges and presence of loops. We present tables of how many
graphs exist in these categories.
\end{abstract}

\maketitle 

\section{Classifications} \label{sec.tag}

A finite graph on $V$ vertices with $E$ edges may
be classified by some properties, which it either does have
or does not:
\begin{itemize}
\item
Each edge in a directed graph has one of two orientations. Edges
in unoriented graphs do not have orientations. 
We reserve the
tag \texttt{d} 
for the directed 
and \texttt{-d} 
for the undirected graphs.
The adjacency matrices of undirected graphs are symmetric.
\item
Graphs may have at least one
loop (loops are defined as edges that start and end at the same vertex), or may be
loopless. The adjacency matrices of loopless graphs have zero trace.
We reserve the
tag \texttt{l} for the graphs with at least one loop
and the tag \texttt{-l} for the loopless graphs.
\item
A multiedge is a collection of two ore more edges having
identical endpoints \cite[D7]{Gross}. This implies that in
a directed graph two edges of opposite sense do not yet establish
a multiedge.
We reserve the
tag \texttt{m} for the graphs with at least one multiedge
and \texttt{-m} for the others.
\item
A undirected graph is connected if one can walk from any vertex
to any other vertex of the graph along edges.
A directed graph is (weakly) connected if replacing each arc with
an undirected edge (defining the underlying graph
\cite[D24]{Gross}) reduces to a connected undirected graph.
This implies that for the sake of weak connectivity it is not required
that all arcs are traversed along their orientation to walk from
one vertex to the other.
We reserve the
tag \texttt{c} for the graphs which are (weakly) connected
and \texttt{-c} for the others.

A directed graph is strongly connected if one can walk from any
vertex to any other vertex of the graph along edges in the directions
demanded by their orientation.
We reserve the
tag \texttt{C} for the digraphs which are strongly connected
and \texttt{-C} for the others.
\item
An isolated vertex is a vertex with no edge to any other vertex (so
all its edges are loops).
There is a loose relation with connectivity, because an isolated vertex
in a graph with two or more vertices means the graph is disconnected.
(There are disconnected graphs without isolated vertices\ldots
where each component contains at least two vertices.)
We reserve the
tag \texttt{i} for the graphs which have at least one isolated vertex
and \texttt{-i} for the others.
Therefore all graphs with $V=1$ are getting the \texttt{i} tag.
\end{itemize}

There are no graphs with the following combinations of tags
\begin{itemize}

\item \texttt{d-Cci} A directed, weakly connected graph 
with at least one isolated vertex has only this vertex
(because with two or more vertices the graph could not be connected), and therefore
the graph must also be strongly connected. So the \texttt{-C}
contradicts the other tags.

\item \texttt{dC-c} If the directed graph is strongly connected,
it is also weakly connected, so the \texttt{-c} tag contradicts
the \texttt{C} tag.

\end{itemize}

There are some non-interesting cases, which are not
tabulated explicitly:
\begin{itemize}
\item
There is the case with the tags \texttt{dCci}:
A directed strongly-connected graph with isolated vertices has
only one vertex, so a table with these graphs counts at most
1 graph for any number of edges (which all are loops).
\item
Similarly there is the case with the tags \texttt{-dci}:
An undirected connected graph with isolated vertices has
only one vertex, so a table with these graphs counts at most
1 graph for any number of edges (where all edges are loops).
\end{itemize}

There are many other characterizations of graphs concerning
cycles, paths, diameters, transitivity and so on which are
not dealt with here.

\clearpage
\section{Statistics}
Tables \ref{tabU.d-C-c-i-m-l}--\ref{tabL.dCc-iml}
collect the statistics of directed graphs;
Tables \ref{tabU.-d-c-i-m-l}--\ref{tabL.-d-ciml}
collect the statistics of unndirected graphs.
Rows and columns are sorted along the number $E$ of edges
and along the number $V$ of vertices.
There are always successive tables referring
to the unlabeled graphs and referring to the
vertex-labeled graphs. (The latter count is obtained
by weighting each unlabeled graph by the number
of distinct adjancency matrices that are created
by row-column permutations of the adjacency matrix.
This weight is a divisor of the order of the permutation
group on $V$ elements \cite[Thm 15.2]{Harary}.)

\subsection{Directed Graphs}
\begin{table}
\begin{tabular}{r|rrrrrrrrrrrrrrrrrr}
$E\backslash V$ & 1 & 2 & 3 & 4 & 5 & 6 & 7 &8 &9 & 10\\
\hline
0& 0& 0& 0& 0& 0& 0& 0& 0& 0& 0& \\
1& 0& 0& 0& 0& 0& 0& 0& 0& 0& 0& \\
2& 0& 0& 0& 1& 0& 0& 0& 0& 0& 0& \\
3& 0& 0& 0& 1& 3& 1& 0& 0& 0& 0& \\
4& 0& 0& 0& 1& 7& 15& 3& 1& 0& 0& \\
5& 0& 0& 0& 0& 8& 43& 58& 15& 3& & \\
6& 0& 0& 0& 0& 5& 82& 244& 257& 68& & \\
7& 0& 0& 0& 0& 2& 103& 674& & & & \\
8& 0& 0& 0& 0& 1& 102& & & & & \\
9& 0& 0& 0& 0& 0& & & & & & \\
10& 0& 0& 0& 0& & & & & & & \\
\hline
\end{tabular}
\caption{\texttt{d-C-c-i-m-l} unlabeled}
\label{tabU.d-C-c-i-m-l}
\end{table}

\begin{table}
\begin{tabular}{r|rrrrrrrrrrrrrrrrrr}
$E\backslash V$ & 1 & 2 & 3 & 4 & 5 & 6 & 7 &8 &9 & 10\\
\hline
0& 0& 0& 0& 0& 0& 0& 0& 0& 0& 0& \\
1& 0& 0& 0& 0& 0& 0& 0& 0& 0& 0& \\
2& 0& 0& 0& 12& 0& 0& 0& 0& 0& 0& \\
3& 0& 0& 0& 12& 240& 120& 0& 0& 0& 0& \\
4& 0& 0& 0& 3& 520& 5460& 5040& 1680& 0& 0& \\
5& 0& 0& 0& 0& 500& 19770& 151200& 191520& 120960& & \\
6& 0& 0& 0& 0& 270& 37135& 795368& 5021912& 7761600& & \\
7& 0& 0& 0& 0& 80& 46560& 2359224& & & & \\
8& 0& 0& 0& 0& 10& 42450& & & & & \\
9& 0& 0& 0& 0& 0& & & & & & \\
10& 0& 0& 0& 0& & & & & & & \\
\hline
\end{tabular}
\caption{\texttt{d-C-c-i-m-l} vertex-labeled}
\label{tabL.d-C-c-i-m-l}
\end{table}

\begin{table}
\begin{tabular}{r|rrrrrrrrrrrrrrrrrr}
$E\backslash V$ & 1 & 2 & 3 & 4 & 5 & 6 & 7 &8 &9 & 10\\
\hline
0& 0& 0& 0& 0& 0& 0& 0& 0& 0& 0& \\
1& 0& 1& 0& 0& 0& 0& 0& 0& 0& 0& \\
2& 0& 0& 3& 0& 0& 0& 0& 0& 0& 0& \\
3& 0& 0& 3& 8& 0& 0& 0& 0& 0& 0& \\
4& 0& 0& 2& 21& 27& 0& 0& 0& 0& 0& \\
5& 0& 0& 0& 33& 107& 91& 0& 0& 0& & \\
6& 0& 0& 0& 31& 319& 581& 350& 0& 0& & \\
7& 0& 0& 0& 16& 609& 2422& 3023& & & & \\
8& 0& 0& 0& 5& 887& 7529& & & & & \\
9& 0& 0& 0& 2& 912& & & & & & \\
10& 0& 0& 0& 0& & & & & & & \\
\hline
\end{tabular}
\caption{\texttt{d-Cc-i-m-l} unlabeled}
\label{tabU.d-Cc-i-m-l}
\end{table}

\begin{table}
\begin{tabular}{r|rrrrrrrrrrrrrrrrrr}
$E\backslash V$ & 1 & 2 & 3 & 4 & 5 & 6 & 7 &8 &9 & 10\\
\hline
0& 0& 0& 0& 0& 0& 0& 0& 0& 0& 0& \\
1& 0& 2& 0& 0& 0& 0& 0& 0& 0& 0& \\
2& 0& 0& 12& 0& 0& 0& 0& 0& 0& 0& \\
3& 0& 0& 18& 128& 0& 0& 0& 0& 0& 0& \\
4& 0& 0& 6& 426& 2000& 0& 0& 0& 0& 0& \\
5& 0& 0& 0& 684& 11080& 41472& 0& 0& 0& & \\
6& 0& 0& 0& 604& 33160& 337800& 1075648& 0& 0& & \\
7& 0& 0& 0& 300& 67040& 1529520& 11967984& & & & \\
8& 0& 0& 0& 78& 96610& 4954230& & & & & \\
9& 0& 0& 0& 8& 101580& & & & & & \\
10& 0& 0& 0& 0& & & & & & & \\
\hline
\end{tabular}
\caption{\texttt{d-Cc-i-m-l} vertex-labeled}
\label{tabL.d-Cc-i-m-l}
\end{table}

\begin{table}
\begin{tabular}{r|rrrrrrrrrrrrrrrrrr}
$E\backslash V$ & 1 & 2 & 3 & 4 & 5 & 6 & 7 &8 &9 & 10\\
\hline
0& 0& 1& 1& 1& 1& 1& 1& 1& 1& 1& \\
1& 0& 0& 1& 1& 1& 1& 1& 1& 1& 1& \\
2& 0& 0& 1& 4& 5& 5& 5& 5& 5& 5& \\
3& 0& 0& 0& 4& 13& 16& 17& 17& 17& 17& \\
4& 0& 0& 0& 4& 27& 61& 76& 79& 80& 80& \\
5& 0& 0& 0& 1& 38& 154& 288& 346& 361& & \\
6& 0& 0& 0& 1& 48& 379& 1043& 1637& 1894& & \\
7& 0& 0& 0& 0& 38& 707& 3242& & & & \\
8& 0& 0& 0& 0& 27& 1155& & & & & \\
9& 0& 0& 0& 0& 13& & & & & & \\
10& 0& 0& 0& 0& & & & & & & \\
\hline
\end{tabular}
\caption{\texttt{d-C-ci-m-l} unlabeled}
\label{tabU.d-C-ci-m-l}
\end{table}

\begin{table}
\begin{tabular}{r|rrrrrrrrrrrrrrrrrr}
$E\backslash V$ & 1 & 2 & 3 & 4 & 5 & 6 & 7 &8 &9 & 10\\
\hline
0& 0& 1& 1& 1& 1& 1& 1& 1& 1& 1& \\
1& 0& 0& 6& 12& 20& 30& 42& 56& 72& 90& \\
2& 0& 0& 3& 54& 190& 435& 861& 1540& 2556& 4005& \\
3& 0& 0& 0& 80& 900& 3940& 11480& 27720& 59640& 117480& \\
4& 0& 0& 0& 60& 2325& 21945& 106890& 365610& 1028790& 2555190& \\
5& 0& 0& 0& 24& 3900& 81264& 699468& 3628296& 13870584& & \\
6& 0& 0& 0& 4& 4610& 218720& 3374770& 27446524& 148477308& & \\
7& 0& 0& 0& 0& 3960& 453240& 12650400& & & & \\
8& 0& 0& 0& 0& 2475& 748395& & & & & \\
9& 0& 0& 0& 0& 1100& & & & & & \\
10& 0& 0& 0& 0& & & & & & & \\
\hline
\end{tabular}
\caption{\texttt{d-C-ci-m-l} vertex-labeled}
\label{tabL.d-C-ci-m-l}
\end{table}

\begin{table}
\begin{tabular}{r|rrrrrrrrrrrrrrrrrr}
$E\backslash V$ & 1 & 2 & 3 & 4 & 5 & 6 & 7 &8 &9 & 10\\
\hline
0& 0& 0& 0& 0& 0& 0& 0& 0& 0& 0& \\
1& 0& 0& 0& 0& 0& 0& 0& 0& 0& 0& \\
2& 0& 0& 0& 0& 0& 0& 0& 0& 0& 0& \\
3& 0& 0& 0& 1& 0& 0& 0& 0& 0& 0& \\
4& 0& 0& 0& 4& 7& 1& 0& 0& 0& 0& \\
5& 0& 0& 0& 7& 35& 42& 7& 1& 0& & \\
6& 0& 0& 0& 12& 101& 271& 234& 48& 7& & \\
7& 0& 0& 0& 16& 230& 1057& 1848& & & & \\
8& 0& 0& 0& 24& 462& 3285& & & & & \\
9& 0& 0& 0& 30& 855& & & & & & \\
10& 0& 0& 0& 41& & & & & & & \\
\hline
\end{tabular}
\caption{\texttt{d-C-c-im-l} unlabeled}
\label{tabU.d-C-c-im-l}
\end{table}

\begin{table}
\begin{tabular}{r|rrrrrrrrrrrrrrrrrr}
$E\backslash V$ & 1 & 2 & 3 & 4 & 5 & 6 & 7 &8 &9 & 10\\
\hline
0& 0& 0& 0& 0& 0& 0& 0& 0& 0& 0& \\
1& 0& 0& 0& 0& 0& 0& 0& 0& 0& 0& \\
2& 0& 0& 0& 0& 0& 0& 0& 0& 0& 0& \\
3& 0& 0& 0& 24& 0& 0& 0& 0& 0& 0& \\
4& 0& 0& 0& 72& 720& 360& 0& 0& 0& 0& \\
5& 0& 0& 0& 132& 3520& 22560& 20160& 6720& 0& & \\
6& 0& 0& 0& 210& 10100& 154650& 806400& 974400& 604800& & \\
7& 0& 0& 0& 312& 23120& 630360& 7141008& & & & \\
8& 0& 0& 0& 441& 46970& 1991325& & & & & \\
9& 0& 0& 0& 600& 88280& & & & & & \\
10& 0& 0& 0& 792& & & & & & & \\
\hline
\end{tabular}
\caption{\texttt{d-C-c-im-l} vertex-labeled}
\label{tabL.d-C-c-im-l}
\end{table}

\clearpage
\begin{table}
\begin{tabular}{r|rrrrrrrrrrrrrrrrrr}
$E\backslash V$ & 1 & 2 & 3 & 4 & 5 & 6 & 7 &8 &9 & 10\\
\hline
0& 0& 0& 0& 0& 0& 0& 0& 0& 0& 0& \\
1& 0& 0& 0& 0& 0& 0& 0& 0& 0& 0& \\
2& 0& 1& 0& 0& 0& 0& 0& 0& 0& 0& \\
3& 0& 1& 4& 0& 0& 0& 0& 0& 0& 0& \\
4& 0& 1& 16& 18& 0& 0& 0& 0& 0& 0& \\
5& 0& 1& 30& 109& 80& 0& 0& 0& 0& & \\
6& 0& 1& 53& 391& 694& 367& 0& 0& 0& & \\
7& 0& 1& 77& 1042& 3574& 4207& 1708& & & & \\
8& 0& 1& 116& 2402& 14093& 29082& & & & & \\
9& 0& 1& 156& 5001& 46144& & & & & & \\
10& 0& 1& 215& 9737& & & & & & & \\
\hline
\end{tabular}
\caption{\texttt{d-Cc-im-l} unlabeled}
\label{tabU.d-Cc-im-l}
\end{table}

\begin{table}
\begin{tabular}{r|rrrrrrrrrrrrrrrrrr}
$E\backslash V$ & 1 & 2 & 3 & 4 & 5 & 6 & 7 &8 &9 & 10\\
\hline
0& 0& 0& 0& 0& 0& 0& 0& 0& 0& 0& \\
1& 0& 0& 0& 0& 0& 0& 0& 0& 0& 0& \\
2& 0& 2& 0& 0& 0& 0& 0& 0& 0& 0& \\
3& 0& 2& 24& 0& 0& 0& 0& 0& 0& 0& \\
4& 0& 2& 90& 384& 0& 0& 0& 0& 0& 0& \\
5& 0& 2& 180& 2472& 8000& 0& 0& 0& 0& & \\
6& 0& 2& 300& 8960& 75400& 207360& 0& 0& 0& & \\
7& 0& 2& 462& 24324& 405160& 2648880& 6453888& & & & \\
8& 0& 2& 672& 56322& 1623440& 19251960& & & & & \\
9& 0& 2& 936& 118168& 5394560& & & & & & \\
10& 0& 2& 1260& 230760& & & & & & & \\
\hline
\end{tabular}
\caption{\texttt{d-Cc-im-l} vertex-labeled}
\label{tabL.d-Cc-im-l}
\end{table}

\begin{table}
\begin{tabular}{r|rrrrrrrrrrrrrrrrrr}
$E\backslash V$ & 1 & 2 & 3 & 4 & 5 & 6 & 7 &8 &9 & 10\\
\hline
0& 0& 0& 0& 0& 0& 0& 0& 0& 0& 0& \\
1& 0& 0& 0& 0& 0& 0& 0& 0& 0& 0& \\
2& 0& 0& 1& 1& 1& 1& 1& 1& 1& 1& \\
3& 0& 0& 2& 6& 7& 7& 7& 7& 7& 7& \\
4& 0& 0& 3& 20& 42& 49& 50& 50& 50& 50& \\
5& 0& 0& 3& 41& 158& 273& 315& 322& 323& & \\
6& 0& 0& 4& 82& 506& 1302& 1940& 2174& 2222& & \\
7& 0& 0& 4& 132& 1330& 5174& 10439& & & & \\
8& 0& 0& 5& 222& 3213& 18293& & & & & \\
9& 0& 0& 5& 335& 7097& & & & & & \\
10& 0& 0& 6& 511& & & & & & & \\
\hline
\end{tabular}
\caption{\texttt{d-C-cim-l} unlabeled}
\label{tabU.d-C-cim-l}
\end{table}

\begin{table}
\begin{tabular}{r|rrrrrrrrrrrrrrrrrr}
$E\backslash V$ & 1 & 2 & 3 & 4 & 5 & 6 & 7 &8 &9 \\
\hline
0& 0& 0& 0& 0& 0& 0& 0& 0& 0\\
1& 0& 0& 0& 0& 0& 0& 0& 0& 0\\
2& 0& 0& 6& 12& 20& 30& 42& 56& 72\\
3& 0& 0& 12& 120& 400& 900& 1764& 3136& 5184\\
4& 0& 0& 15& 414& 3290& 13155& 37065& 87836& 186660\\
5& 0& 0& 18& 948& 15480& 113190& 499926& 1634976& 4483296& & \\
6& 0& 0& 21& 1802& 52720& 667375& 4685387& 22082536& 80250072& & \\
7& 0& 0& 24& 3120& 147320& 3031920& 33055848& & & & \\
8& 0& 0& 27& 5094& 362655& 11463930& & & & & \\
9& 0& 0& 30& 7948& 818780& & & & & & \\
10& 0& 0& 33& 11946& & & & & & & \\
\hline
\end{tabular}
\caption{\texttt{d-C-cim-l} vertex-labeled}
\label{tabL.d-C-cim-l}
\end{table}

\begin{table}
\begin{tabular}{r|rrrrrrrrrrrrrrrrrr}
$E\backslash V$ & 1 & 2 & 3 & 4 & 5 & 6 & 7 &8 &9 & 10\\
\hline
0& 0& 0& 0& 0& 0& 0& 0& 0& 0& 0& \\
1& 0& 0& 0& 0& 0& 0& 0& 0& 0& 0& \\
2& 0& 0& 0& 0& 0& 0& 0& 0& 0& 0& \\
3& 0& 0& 0& 2& 0& 0& 0& 0& 0& 0& \\
4& 0& 0& 0& 7& 13& 2& 0& 0& 0& 0& \\
5& 0& 0& 0& 7& 52& 70& 13& 2& 0& & \\
6& 0& 0& 0& 6& 106& 373& 362& 82& 13& & \\
7& 0& 0& 0& 2& 137& 1092& 2392& & & & \\
8& 0& 0& 0& 1& 125& 2262& & & & & \\
9& 0& 0& 0& 0& 83& & & & & & \\
10& 0& 0& 0& 0& & & & & & & \\
\hline
\end{tabular}
\caption{\texttt{d-C-c-i-ml} unlabeled}
\label{tabU.d-C-c-i-ml}
\end{table}

\begin{table}
\begin{tabular}{r|rrrrrrrrrrrrrrrrrr}
$E\backslash V$ & 1 & 2 & 3 & 4 & 5 & 6 & 7 &8 &9 & 10\\
\hline
0& 0& 0& 0& 0& 0& 0& 0& 0& 0& 0& \\
1& 0& 0& 0& 0& 0& 0& 0& 0& 0& 0& \\
2& 0& 0& 0& 0& 0& 0& 0& 0& 0& 0& \\
3& 0& 0& 0& 48& 0& 0& 0& 0& 0& 0& \\
4& 0& 0& 0& 120& 1200& 720& 0& 0& 0& 0& \\
5& 0& 0& 0& 132& 5000& 34560& 35280& 13440& 0& & \\
6& 0& 0& 0& 78& 10100& 202920& 1164240& 1579200& 1088640& & \\
7& 0& 0& 0& 24& 12750& 630360& 8919176& & & & \\
8& 0& 0& 0& 3& 10940& 1314405& & & & & \\
9& 0& 0& 0& 0& 6570& & & & & & \\
10& 0& 0& 0& 0& & & & & & & \\
\hline
\end{tabular}
\caption{\texttt{d-C-c-i-ml} vertex-labeled}
\label{tabL.d-C-c-i-ml}
\end{table}

\begin{table}
\begin{tabular}{r|rrrrrrrrrrrrrrrrrr}
$E\backslash V$ & 1 & 2 & 3 & 4 & 5 & 6 & 7 &8 &9 & 10\\
\hline
0& 0& 0& 0& 0& 0& 0& 0& 0& 0& 0& \\
1& 0& 0& 0& 0& 0& 0& 0& 0& 0& 0& \\
2& 0& 2& 0& 0& 0& 0& 0& 0& 0& 0& \\
3& 0& 1& 7& 0& 0& 0& 0& 0& 0& 0& \\
4& 0& 0& 16& 26& 0& 0& 0& 0& 0& 0& \\
5& 0& 0& 16& 111& 107& 0& 0& 0& 0& & \\
6& 0& 0& 7& 262& 702& 458& 0& 0& 0& & \\
7& 0& 0& 2& 372& 2663& 4251& 2058& & & & \\
8& 0& 0& 0& 361& 6936& 22925& & & & & \\
9& 0& 0& 0& 240& 13442& & & & & & \\
10& 0& 0& 0& 115& & & & & & & \\
\hline
\end{tabular}
\caption{\texttt{d-Cc-i-ml} unlabeled}
\label{tabU.d-Cc-i-ml}
\end{table}

\begin{table}
\begin{tabular}{r|rrrrrrrrrrrrrrrrrr}
$E\backslash V$ & 1 & 2 & 3 & 4 & 5 & 6 & 7 &8 &9 & 10\\
\hline
0& 0& 0& 0& 0& 0& 0& 0& 0& 0& 0& \\
1& 0& 0& 0& 0& 0& 0& 0& 0& 0& 0& \\
2& 0& 4& 0& 0& 0& 0& 0& 0& 0& 0& \\
3& 0& 2& 36& 0& 0& 0& 0& 0& 0& 0& \\
4& 0& 0& 90& 512& 0& 0& 0& 0& 0& 0& \\
5& 0& 0& 84& 2472& 10000& 0& 0& 0& 0& & \\
6& 0& 0& 36& 5804& 75400& 248832& 0& 0& 0& & \\
7& 0& 0& 6& 8352& 296600& 2648880& 7529536& & & & \\
8& 0& 0& 0& 7986& 787600& 15073560& & & & & \\
9& 0& 0& 0& 5212& 1542450& & & & & & \\
10& 0& 0& 0& 2304& & & & & & & \\
\hline
\end{tabular}
\caption{\texttt{d-Cc-i-ml} vertex-labeled}
\label{tabL.d-Cc-i-ml}
\end{table}

\clearpage
\begin{table}
\begin{tabular}{r|rrrrrrrrrrrrrrrrrr}
$E\backslash V$ & 1 & 2 & 3 & 4 & 5 & 6 & 7 &8 &9 & 10\\
\hline
0& 0& 0& 0& 0& 0& 0& 0& 0& 0& 0& \\
1& 0& 1& 1& 1& 1& 1& 1& 1& 1& 1& \\
2& 0& 1& 4& 4& 4& 4& 4& 4& 4& 4& \\
3& 0& 0& 6& 17& 20& 20& 20& 20& 20& 20& \\
4& 0& 0& 3& 35& 83& 100& 103& 103& 103& 103& \\
5& 0& 0& 1& 46& 236& 457& 548& 565& 568& & \\
6& 0& 0& 0& 40& 504& 1659& 2756& 3210& 3313& & \\
7& 0& 0& 0& 25& 833& 4986& 12171& & & & \\
8& 0& 0& 0& 10& 1064& 12330& & & & & \\
9& 0& 0& 0& 3& 1084& & & & & & \\
10& 0& 0& 0& 1& & & & & & & \\
\hline
\end{tabular}
\caption{\texttt{d-C-ci-ml} unlabeled}
\label{tabU.d-C-ci-ml}
\end{table}

\begin{table}
\begin{tabular}{r|rrrrrrrrrrrrrrrrrr}
$E\backslash V$ & 1 & 2 & 3 & 4 & 5 & 6 & 7 &8 &9 \\
\hline
0& 0& 0& 0& 0& 0& 0& 0& 0& 0\\
1& 0& 2& 3& 4& 5& 6& 7& 8& 9\\
2& 0& 1& 21& 54& 110& 195& 315& 476& 684\\
3& 0& 0& 28& 292& 1160& 3080& 6944& 13944& 25680\\
4& 0& 0& 15& 693& 6605& 30780& 99946& 268086& 634950\\
5& 0& 0& 3& 948& 22626& 199926& 1020936& 3791256& 11630052& & \\
6& 0& 0& 0& 830& 52720& 902265& 7573790& 40926732& 167212668& & \\
7& 0& 0& 0& 480& 89990& 3031920& 42473544& & & & \\
8& 0& 0& 0& 180& 117425& 7978770& & & & & \\
9& 0& 0& 0& 40& 119900& & & & & & \\
10& 0& 0& 0& 4& & & & & & & \\
\hline
\end{tabular}
\caption{\texttt{d-C-ci-ml} vertex-labeled}
\label{tabL.d-C-ci-ml}
\end{table}

\begin{table}
\begin{tabular}{r|rrrrrrrrrrrrrrrrrr}
$E\backslash V$ & 1 & 2 & 3 & 4 & 5 & 6 & 7 &8 &9 & 10\\
\hline
0& 0& 0& 0& 0& 0& 0& 0& 0& 0& 0& \\
1& 0& 0& 0& 0& 0& 0& 0& 0& 0& 0& \\
2& 0& 0& 0& 0& 0& 0& 0& 0& 0& 0& \\
3& 0& 0& 0& 0& 0& 0& 0& 0& 0& 0& \\
4& 0& 0& 0& 6& 0& 0& 0& 0& 0& 0& \\
5& 0& 0& 0& 34& 46& 6& 0& 0& 0& & \\
6& 0& 0& 0& 107& 347& 314& 52& 6& 0& & \\
7& 0& 0& 0& 250& 1473& 2869& 1995& & & & \\
8& 0& 0& 0& 527& 4731& 15676& & & & & \\
9& 0& 0& 0& 994& 12883& & & & & & \\
10& 0& 0& 0& 1797& & & & & & & \\
\hline
\end{tabular}
\caption{\texttt{d-C-c-iml} unlabeled}
\label{tabU.d-C-c-iml}
\end{table}

\begin{table}
\begin{tabular}{r|rrrrrrrrrrrrrrrrrr}
$E\backslash V$ & 1 & 2 & 3 & 4 & 5 & 6 & 7 &8 &9 & 10\\
\hline
0& 0& 0& 0& 0& 0& 0& 0& 0& 0& 0& \\
1& 0& 0& 0& 0& 0& 0& 0& 0& 0& 0& \\
2& 0& 0& 0& 0& 0& 0& 0& 0& 0& 0& \\
3& 0& 0& 0& 0& 0& 0& 0& 0& 0& 0& \\
4& 0& 0& 0& 144& 0& 0& 0& 0& 0& 0& \\
5& 0& 0& 0& 768& 4800& 2880& 0& 0& 0& & \\
6& 0& 0& 0& 2340& 37000& 180000& 176400& 67200& 0& & \\
7& 0& 0& 0& 5568& 159600& 1750320& 7514640& & & & \\
8& 0& 0& 0& 11634& 518350& 9908640& & & & & \\
9& 0& 0& 0& 22368& 1427320& & & & & & \\
10& 0& 0& 0& 40392& & & & & & & \\
\hline
\end{tabular}
\caption{\texttt{d-C-c-iml} vertex-labeled}
\label{tabL.d-C-c-iml}
\end{table}

\begin{table}
\begin{tabular}{r|rrrrrrrrrrrrrrrrrr}
$E\backslash V$ & 1 & 2 & 3 & 4 & 5 & 6 & 7 &8 &9 & 10\\
\hline
0& 0& 0& 0& 0& 0& 0& 0& 0& 0& 0& \\
1& 0& 0& 0& 0& 0& 0& 0& 0& 0& 0& \\
2& 0& 0& 0& 0& 0& 0& 0& 0& 0& 0& \\
3& 0& 4& 0& 0& 0& 0& 0& 0& 0& 0& \\
4& 0& 9& 19& 0& 0& 0& 0& 0& 0& 0& \\
5& 0& 14& 98& 94& 0& 0& 0& 0& 0& & \\
6& 0& 20& 286& 761& 479& 0& 0& 0& 0& & \\
7& 0& 27& 645& 3522& 5398& 2480& 0& & & & \\
8& 0& 35& 1290& 12111& 34960& 36619& & & & & \\
9& 0& 44& 2372& 34847& 167682& & & & & & \\
10& 0& 54& 4110& 89361& & & & & & & \\
\hline
\end{tabular}
\caption{\texttt{d-Cc-iml} unlabeled}
\label{tabU.d-Cc-iml}
\end{table}

\begin{table}
\begin{tabular}{r|rrrrrrrrrrrrrrrrrr}
$E\backslash V$ & 1 & 2 & 3 & 4 & 5 & 6 & 7 &8 &9 & 10\\
\hline
0& 0& 0& 0& 0& 0& 0& 0& 0& 0& 0& \\
1& 0& 0& 0& 0& 0& 0& 0& 0& 0& 0& \\
2& 0& 0& 0& 0& 0& 0& 0& 0& 0& 0& \\
3& 0& 8& 0& 0& 0& 0& 0& 0& 0& 0& \\
4& 0& 18& 108& 0& 0& 0& 0& 0& 0& 0& \\
5& 0& 28& 576& 2048& 0& 0& 0& 0& 0& & \\
6& 0& 40& 1680& 17480& 50000& 0& 0& 0& 0& & \\
7& 0& 54& 3816& 82144& 602400& 1492992& 0& & & & \\
8& 0& 70& 7644& 284788& 4009600& 23767632& & & & & \\
9& 0& 88& 14112& 824480& 19528000& & & & & & \\
10& 0& 108& 24480& 2121756& & & & & & & \\
\hline
\end{tabular}
\caption{\texttt{d-Cc-iml} vertex-labeled}
\label{tabL.d-Cc-iml}
\end{table}

\begin{table}
\begin{tabular}{r|rrrrrrrrrrrrrrrrrr}
$E\backslash V$ & 1 & 2 & 3 & 4 & 5 & 6 & 7 &8 &9 & 10\\
\hline
0& 0& 0& 0& 0& 0& 0& 0& 0& 0& 0& \\
1& 0& 0& 0& 0& 0& 0& 0& 0& 0& 0& \\
2& 0& 1& 1& 1& 1& 1& 1& 1& 1& 1& \\
3& 0& 2& 8& 8& 8& 8& 8& 8& 8& 8& \\
4& 0& 3& 27& 55& 63& 63& 63& 63& 63& 63& \\
5& 0& 3& 55& 224& 402& 460& 468& 468& 468& & \\
6& 0& 4& 97& 671& 1956& 3051& 3444& 3508& 3516& & \\
7& 0& 4& 154& 1661& 7607& 17024& 23868& & & & \\
8& 0& 5& 235& 3670& 25207& 81289& & & & & \\
9& 0& 5& 342& 7505& 74029& & & & & & \\
10& 0& 6& 483& 14483& & & & & & & \\
\hline
\end{tabular}
\caption{\texttt{d-C-ciml} unlabeled}
\label{tabU.d-C-ciml}
\end{table}

\begin{table}
\begin{tabular}{r|rrrrrrrrrrrrrrrrrr}
$E\backslash V$ & 1 & 2 & 3 & 4 & 5 & 6 & 7 &8 &9 \\
\hline
0& 0& 0& 0& 0& 0& 0& 0& 0& 0\\
1& 0& 0& 0& 0& 0& 0& 0& 0& 0\\
2& 0& 2& 3& 4& 5& 6& 7& 8& 9\\
3& 0& 4& 45& 112& 225& 396& 637& 960& 1377\\
4& 0& 5& 150& 1042& 3815& 9831& 21784& 43268& 79101\\
5& 0& 6& 315& 4744& 33825& 142386& 442715& 1157920& 2696625& & \\
6& 0& 7& 553& 14864& 191335& 1339211& 6175162& 21778968& 64759989& & \\
7& 0& 8& 894& 37768& 805875& 9075456& 62861757& & & & \\
8& 0& 9& 1368& 84739& 2785270& 48311766& & & & & \\
9& 0& 10& 2005& 175140& 8398505& & & & & & \\
10& 0& 11& 2838& 340402& & & & & & & \\
\hline
\end{tabular}
\caption{\texttt{d-C-ciml} vertex-labeled}
\label{tabL.d-C-ciml}
\end{table}

\clearpage
\begin{table}
\begin{tabular}{r|rrrrrrrrrrrrrrrrrr}
$E\backslash V$ & 1 & 2 & 3 & 4 & 5 & 6 & 7 &8 &9 & 10\\
\hline
0& 0& 0& 0& 0& 0& 0& 0& 0& 0& 0& \\
1& 0& 0& 0& 0& 0& 0& 0& 0& 0& 0& \\
2& 0& 1& 0& 0& 0& 0& 0& 0& 0& 0& \\
3& 0& 0& 1& 0& 0& 0& 0& 0& 0& 0& \\
4& 0& 0& 2& 1& 0& 0& 0& 0& 0& 0& \\
5& 0& 0& 1& 4& 1& 0& 0& 0& 0& & \\
6& 0& 0& 1& 16& 7& 1& 0& 0& 0& & \\
7& 0& 0& 0& 22& 58& 10& 1& & & & \\
8& 0& 0& 0& 22& 240& 165& & & & & \\
9& 0& 0& 0& 11& 565& & & & & & \\
10& 0& 0& 0& 5& & & & & & & \\
\hline
\end{tabular}
\caption{\texttt{dCc-i-m-l} unlabeled}
\label{tabU.dCc-i-m-l}
\end{table}

\begin{table}
\begin{tabular}{r|rrrrrrrrrrrrrrrrrr}
$E\backslash V$ & 1 & 2 & 3 & 4 & 5 & 6 & 7 &8 &9 & 10\\
\hline
0& 0& 0& 0& 0& 0& 0& 0& 0& 0& 0& \\
1& 0& 0& 0& 0& 0& 0& 0& 0& 0& 0& \\
2& 0& 1& 0& 0& 0& 0& 0& 0& 0& 0& \\
3& 0& 0& 2& 0& 0& 0& 0& 0& 0& 0& \\
4& 0& 0& 9& 6& 0& 0& 0& 0& 0& 0& \\
5& 0& 0& 6& 84& 24& 0& 0& 0& 0& & \\
6& 0& 0& 1& 316& 720& 120& 0& 0& 0& & \\
7& 0& 0& 0& 492& 6440& 6480& 720& & & & \\
8& 0& 0& 0& 417& 26875& 107850& & & & & \\
9& 0& 0& 0& 212& 65280& & & & & & \\
10& 0& 0& 0& 66& & & & & & & \\
\hline
\end{tabular}
\caption{\texttt{dCc-i-m-l} vertex-labeled}
\label{tabL.dCc-i-m-l}
\end{table}

\begin{table}
\begin{tabular}{r|rrrrrrrrrrrrrrrrrr}
$E\backslash V$ & 1 & 2 & 3 & 4 & 5 & 6 & 7 &8 &9 & 10\\
\hline
0& 0& 0& 0& 0& 0& 0& 0& 0& 0& 0& \\
1& 0& 0& 0& 0& 0& 0& 0& 0& 0& 0& \\
2& 0& 0& 0& 0& 0& 0& 0& 0& 0& 0& \\
3& 0& 1& 0& 0& 0& 0& 0& 0& 0& 0& \\
4& 0& 2& 1& 0& 0& 0& 0& 0& 0& 0& \\
5& 0& 2& 8& 1& 0& 0& 0& 0& 0& & \\
6& 0& 3& 25& 21& 1& 0& 0& 0& 0& & \\
7& 0& 3& 51& 140& 40& 1& 0& & & & \\
8& 0& 4& 101& 565& 525& 69& & & & & \\
9& 0& 4& 174& 1731& 3719& & & & & & \\
10& 0& 5& 290& 4602& & & & & & & \\
\hline
\end{tabular}
\caption{\texttt{dCc-im-l} unlabeled}
\label{tabU.dCc-im-l}
\end{table}

\begin{table}
\begin{tabular}{r|rrrrrrrrrrrrrrrrrr}
$E\backslash V$ & 1 & 2 & 3 & 4 & 5 & 6 & 7 &8 &9 & 10\\
\hline
0& 0& 0& 0& 0& 0& 0& 0& 0& 0& 0& \\
1& 0& 0& 0& 0& 0& 0& 0& 0& 0& 0& \\
2& 0& 0& 0& 0& 0& 0& 0& 0& 0& 0& \\
3& 0& 2& 0& 0& 0& 0& 0& 0& 0& 0& \\
4& 0& 3& 6& 0& 0& 0& 0& 0& 0& 0& \\
5& 0& 4& 48& 24& 0& 0& 0& 0& 0& & \\
6& 0& 5& 140& 480& 120& 0& 0& 0& 0& & \\
7& 0& 6& 306& 3276& 4680& 720& 0& & & & \\
8& 0& 7& 588& 13230& 61040& 47880& & & & & \\
9& 0& 8& 1036& 41024& 437320& & & & & & \\
10& 0& 9& 1710& 109152& & & & & & & \\
\hline
\end{tabular}
\caption{\texttt{dCc-im-l} vertex-labeled}
\label{tabL.dCc-im-l}
\end{table}

\begin{table}
\begin{tabular}{r|rrrrrrrrrrrrrrrrrr}
$E\backslash V$ & 1 & 2 & 3 & 4 & 5 & 6 & 7 &8 &9 & 10\\
\hline
0& 0& 0& 0& 0& 0& 0& 0& 0& 0& 0& \\
1& 0& 0& 0& 0& 0& 0& 0& 0& 0& 0& \\
2& 0& 0& 0& 0& 0& 0& 0& 0& 0& 0& \\
3& 0& 1& 0& 0& 0& 0& 0& 0& 0& 0& \\
4& 0& 1& 1& 0& 0& 0& 0& 0& 0& 0& \\
5& 0& 0& 6& 1& 0& 0& 0& 0& 0& & \\
6& 0& 0& 9& 17& 1& 0& 0& 0& 0& & \\
7& 0& 0& 6& 78& 34& 1& 0& & & & \\
8& 0& 0& 2& 185& 346& 60& & & & & \\
9& 0& 0& 1& 259& 1775& & & & & & \\
10& 0& 0& 0& 252& & & & & & & \\
\hline
\end{tabular}
\caption{\texttt{dCc-i-ml} unlabeled}
\label{tabU.dCc-i-ml}
\end{table}

\begin{table}
\begin{tabular}{r|rrrrrrrrrrrrrrrrrr}
$E\backslash V$ & 1 & 2 & 3 & 4 & 5 & 6 & 7 &8 &9 & 10\\
\hline
0& 0& 0& 0& 0& 0& 0& 0& 0& 0& 0& \\
1& 0& 0& 0& 0& 0& 0& 0& 0& 0& 0& \\
2& 0& 0& 0& 0& 0& 0& 0& 0& 0& 0& \\
3& 0& 2& 0& 0& 0& 0& 0& 0& 0& 0& \\
4& 0& 1& 6& 0& 0& 0& 0& 0& 0& 0& \\
5& 0& 0& 33& 24& 0& 0& 0& 0& 0& & \\
6& 0& 0& 47& 372& 120& 0& 0& 0& 0& & \\
7& 0& 0& 30& 1792& 3840& 720& 0& & & & \\
8& 0& 0& 9& 4206& 39640& 40680& & & & & \\
9& 0& 0& 1& 5968& 206095& & & & & & \\
10& 0& 0& 0& 5634& & & & & & & \\
\hline
\end{tabular}
\caption{\texttt{dCc-i-ml} vertex-labeled}
\label{tabL.dCc-i-ml}
\end{table}

\begin{table}
\begin{tabular}{r|rrrrrrrrrrrrrrrrrr}
$E\backslash V$ & 1 & 2 & 3 & 4 & 5 & 6 & 7 &8 &9 & 10\\
\hline
0& 0& 0& 0& 0& 0& 0& 0& 0& 0& 0& \\
1& 0& 0& 0& 0& 0& 0& 0& 0& 0& 0& \\
2& 0& 0& 0& 0& 0& 0& 0& 0& 0& 0& \\
3& 0& 0& 0& 0& 0& 0& 0& 0& 0& 0& \\
4& 0& 3& 0& 0& 0& 0& 0& 0& 0& 0& \\
5& 0& 8& 4& 0& 0& 0& 0& 0& 0& & \\
6& 0& 16& 38& 5& 0& 0& 0& 0& 0& & \\
7& 0& 25& 151& 110& 6& 0& 0& & & & \\
8& 0& 40& 431& 898& 250& 7& & & & & \\
9& 0& 56& 1040& 4475& 3665& & & & & & \\
10& 0& 80& 2252& 17039& & & & & & & \\
\hline
\end{tabular}
\caption{\texttt{dCc-iml} unlabeled}
\label{tabU.dCc-iml}
\end{table}

\begin{table}
\begin{tabular}{r|rrrrrrrrrrrrrrrrrr}
$E\backslash V$ & 1 & 2 & 3 & 4 & 5 & 6 & 7 &8 &9 & 10\\
\hline
0& 0& 0& 0& 0& 0& 0& 0& 0& 0& 0& \\
1& 0& 0& 0& 0& 0& 0& 0& 0& 0& 0& \\
2& 0& 0& 0& 0& 0& 0& 0& 0& 0& 0& \\
3& 0& 0& 0& 0& 0& 0& 0& 0& 0& 0& \\
4& 0& 6& 0& 0& 0& 0& 0& 0& 0& 0& \\
5& 0& 16& 24& 0& 0& 0& 0& 0& 0& & \\
6& 0& 30& 225& 120& 0& 0& 0& 0& 0& & \\
7& 0& 50& 897& 2592& 720& 0& 0& & & & \\
8& 0& 77& 2562& 21196& 29400& 5040& & & & & \\
9& 0& 112& 6190& 106336& 431360& & & & & & \\
10& 0& 156& 13437& 405552& & & & & & & \\
\hline
\end{tabular}
\caption{\texttt{dCc-iml} vertex-labeled}
\label{tabL.dCc-iml}
\end{table}

\clearpage
\subsection{Undirected Graphs}
\begin{table}
\begin{tabular}{r|rrrrrrrrrrrrrrrrrr}
$E\backslash V$ & 1 & 2 & 3 & 4 & 5 & 6 & 7 &8 &9 & 10\\
\hline
0& 0& 0& 0& 0& 0& 0& 0& 0& 0& 0& \\
1& 0& 0& 0& 0& 0& 0& 0& 0& 0& 0& \\
2& 0& 0& 0& 1& 0& 0& 0& 0& 0& 0& \\
3& 0& 0& 0& 0& 1& 1& 0& 0& 0& 0& \\
4& 0& 0& 0& 0& 1& 3& 1& 1& 0& 0& \\
5& 0& 0& 0& 0& 0& 3& 6& 3& 1& & \\
6& 0& 0& 0& 0& 0& 2& 9& 15& 7& & \\
7& 0& 0& 0& 0& 0& 1& 8& & & & \\
8& 0& 0& 0& 0& 0& 0& & & & & \\
9& 0& 0& 0& 0& 0& & & & & & \\
10& 0& 0& 0& 0& & & & & & & \\
\hline
\end{tabular}
\caption{\texttt{-d-c-i-m-l} unlabeled}
\label{tabU.-d-c-i-m-l}
\end{table}

\begin{table}
\begin{tabular}{r|rrrrrrrrrrrrrrrrrr}
$E\backslash V$ & 1 & 2 & 3 & 4 & 5 & 6 & 7 &8 &9 & 10\\
\hline
0& 0& 0& 0& 0& 0& 0& 0& 0& 0& 0& \\
1& 0& 0& 0& 0& 0& 0& 0& 0& 0& 0& \\
2& 0& 0& 0& 3& 0& 0& 0& 0& 0& 0& \\
3& 0& 0& 0& 0& 30& 15& 0& 0& 0& 0& \\
4& 0& 0& 0& 0& 10& 330& 315& 105& 0& 0& \\
5& 0& 0& 0& 0& 0& 285& 4410& 5880& 3780& & \\
6& 0& 0& 0& 0& 0& 100& 6797& 71078& 116550& & \\
7& 0& 0& 0& 0& 0& 15& 5460& & & & \\
8& 0& 0& 0& 0& 0& 0& & & & & \\
9& 0& 0& 0& 0& 0& & & & & & \\
10& 0& 0& 0& 0& & & & & & & \\
\hline
\end{tabular}
\caption{\texttt{-d-c-i-m-l} vertex-labeled}
\label{tabL.-d-c-i-m-l}
\end{table}

\begin{table}
\begin{tabular}{r|rrrrrrrrrrrrrrrrrr}
$E\backslash V$ & 1 & 2 & 3 & 4 & 5 & 6 & 7 &8 &9 & 10\\
\hline
0& 0& 0& 0& 0& 0& 0& 0& 0& 0& 0& \\
1& 0& 1& 0& 0& 0& 0& 0& 0& 0& 0& \\
2& 0& 0& 1& 0& 0& 0& 0& 0& 0& 0& \\
3& 0& 0& 1& 2& 0& 0& 0& 0& 0& 0& \\
4& 0& 0& 0& 2& 3& 0& 0& 0& 0& 0& \\
5& 0& 0& 0& 1& 5& 6& 0& 0& 0& & \\
6& 0& 0& 0& 1& 5& 13& 11& 0& 0& & \\
7& 0& 0& 0& 0& 4& 19& 33& & & & \\
8& 0& 0& 0& 0& 2& 22& & & & & \\
9& 0& 0& 0& 0& 1& & & & & & \\
10& 0& 0& 0& 0& & & & & & & \\
\hline
\end{tabular}
\caption{\texttt{-dc-i-m-l} unlabeled}
\label{tabU.-dc-i-m-l}
\end{table}

\begin{table}
\begin{tabular}{r|rrrrrrrrrrrrrrrrrr}
$E\backslash V$ & 1 & 2 & 3 & 4 & 5 & 6 & 7 &8 &9 & 10\\
\hline
0& 0& 0& 0& 0& 0& 0& 0& 0& 0& 0& \\
1& 0& 1& 0& 0& 0& 0& 0& 0& 0& 0& \\
2& 0& 0& 3& 0& 0& 0& 0& 0& 0& 0& \\
3& 0& 0& 1& 16& 0& 0& 0& 0& 0& 0& \\
4& 0& 0& 0& 15& 125& 0& 0& 0& 0& 0& \\
5& 0& 0& 0& 6& 222& 1296& 0& 0& 0& & \\
6& 0& 0& 0& 1& 205& 3660& 16807& 0& 0& & \\
7& 0& 0& 0& 0& 120& 5700& 68295& & & & \\
8& 0& 0& 0& 0& 45& 6165& & & & & \\
9& 0& 0& 0& 0& 10& & & & & & \\
10& 0& 0& 0& 0& & & & & & & \\
\hline
\end{tabular}
\caption{\texttt{-dc-i-m-l} vertex-labeled}
\label{tabL.-dc-i-m-l}
\end{table}

\begin{table}
\begin{tabular}{r|rrrrrrrrrrrrrrrrrr}
$E\backslash V$ & 1 & 2 & 3 & 4 & 5 & 6 & 7 &8 &9 & 10\\
\hline
0& 0& 1& 1& 1& 1& 1& 1& 1& 1& 1& \\
1& 0& 0& 1& 1& 1& 1& 1& 1& 1& 1& \\
2& 0& 0& 0& 1& 2& 2& 2& 2& 2& 2& \\
3& 0& 0& 0& 1& 3& 4& 5& 5& 5& 5& \\
4& 0& 0& 0& 0& 2& 6& 9& 10& 11& 11& \\
5& 0& 0& 0& 0& 1& 6& 15& 21& 24& & \\
6& 0& 0& 0& 0& 1& 6& 21& 41& 56& & \\
7& 0& 0& 0& 0& 0& 4& 24& & & & \\
8& 0& 0& 0& 0& 0& 2& & & & & \\
9& 0& 0& 0& 0& 0& & & & & & \\
10& 0& 0& 0& 0& & & & & & & \\
\hline
\end{tabular}
\caption{\texttt{-d-ci-m-l} unlabeled}
\label{tabU.-d-ci-m-l}
\end{table}

\begin{table}
\begin{tabular}{r|rrrrrrrrrrrrrrrrrr}
$E\backslash V$ & 1 & 2 & 3 & 4 & 5 & 6 & 7 &8 &9 & 10\\
\hline
0& 0& 1& 1& 1& 1& 1& 1& 1& 1& 1& \\
1& 0& 0& 3& 6& 10& 15& 21& 28& 36& 45& \\
2& 0& 0& 0& 12& 45& 105& 210& 378& 630& 990& \\
3& 0& 0& 0& 4& 90& 440& 1330& 3276& 7140& 14190& \\
4& 0& 0& 0& 0& 75& 1035& 5670& 20370& 58905& 148995& \\
5& 0& 0& 0& 0& 30& 1422& 15939& 92400& 373212& & \\
6& 0& 0& 0& 0& 5& 1245& 30660& 305662& 1831242& & \\
7& 0& 0& 0& 0& 0& 720& 42525& & & & \\
8& 0& 0& 0& 0& 0& 270& & & & & \\
9& 0& 0& 0& 0& 0& & & & & & \\
10& 0& 0& 0& 0& & & & & & & \\
\hline
\end{tabular}
\caption{\texttt{-d-ci-m-l} vertex-labeled}
\label{tabL.-d-ci-m-l}
\end{table}

\begin{table}
\begin{tabular}{r|rrrrrrrrrrrrrrrrrr}
$E\backslash V$ & 1 & 2 & 3 & 4 & 5 & 6 & 7 &8 &9 & 10\\
\hline
0& 0& 0& 0& 0& 0& 0& 0& 0& 0& 0& \\
1& 0& 0& 0& 0& 0& 0& 0& 0& 0& 0& \\
2& 0& 0& 0& 0& 0& 0& 0& 0& 0& 0& \\
3& 0& 0& 0& 1& 0& 0& 0& 0& 0& 0& \\
4& 0& 0& 0& 2& 2& 1& 0& 0& 0& 0& \\
5& 0& 0& 0& 2& 6& 8& 2& 1& 0& & \\
6& 0& 0& 0& 3& 10& 25& 21& 9& 2& & \\
7& 0& 0& 0& 3& 16& 53& 80& & & & \\
8& 0& 0& 0& 4& 23& 102& & & & & \\
9& 0& 0& 0& 4& 32& & & & & & \\
10& 0& 0& 0& 5& & & & & & & \\
\hline
\end{tabular}
\caption{\texttt{-d-c-im-l} unlabeled}
\label{tabU.-d-c-im-l}
\end{table}

\begin{table}
\begin{tabular}{r|rrrrrrrrrrrrrrrrrr}
$E\backslash V$ & 1 & 2 & 3 & 4 & 5 & 6 & 7 &8 &9 & 10\\
\hline
0& 0& 0& 0& 0& 0& 0& 0& 0& 0& 0& \\
1& 0& 0& 0& 0& 0& 0& 0& 0& 0& 0& \\
2& 0& 0& 0& 0& 0& 0& 0& 0& 0& 0& \\
3& 0& 0& 0& 6& 0& 0& 0& 0& 0& 0& \\
4& 0& 0& 0& 9& 90& 45& 0& 0& 0& 0& \\
5& 0& 0& 0& 12& 220& 1410& 1260& 420& 0& & \\
6& 0& 0& 0& 15& 400& 4875& 25200& 30450& 18900& & \\
7& 0& 0& 0& 18& 650& 11700& 113232& & & & \\
8& 0& 0& 0& 21& 980& 24045& & & & & \\
9& 0& 0& 0& 24& 1400& & & & & & \\
10& 0& 0& 0& 27& & & & & & & \\
\hline
\end{tabular}
\caption{\texttt{-d-c-im-l} vertex-labeled}
\label{tabL.-d-c-im-l}
\end{table}

\clearpage
\begin{table}
\begin{tabular}{r|rrrrrrrrrrrrrrrrrr}
$E\backslash V$ & 1 & 2 & 3 & 4 & 5 & 6 & 7 &8 &9 & 10\\
\hline
0& 0& 0& 0& 0& 0& 0& 0& 0& 0& 0& \\
1& 0& 0& 0& 0& 0& 0& 0& 0& 0& 0& \\
2& 0& 1& 0& 0& 0& 0& 0& 0& 0& 0& \\
3& 0& 1& 1& 0& 0& 0& 0& 0& 0& 0& \\
4& 0& 1& 3& 3& 0& 0& 0& 0& 0& 0& \\
5& 0& 1& 4& 10& 6& 0& 0& 0& 0& & \\
6& 0& 1& 6& 21& 29& 16& 0& 0& 0& & \\
7& 0& 1& 7& 37& 81& 91& 37& & & & \\
8& 0& 1& 9& 61& 191& 326& & & & & \\
9& 0& 1& 11& 95& 395& & & & & & \\
10& 0& 1& 13& 141& & & & & & & \\
\hline
\end{tabular}
\caption{\texttt{-dc-im-l} unlabeled}
\label{tabU.-dc-im-l}
\end{table}

\begin{table}
\begin{tabular}{r|rrrrrrrrrrrrrrrrrr}
$E\backslash V$ & 1 & 2 & 3 & 4 & 5 & 6 & 7 &8 &9 & 10\\
\hline
0& 0& 0& 0& 0& 0& 0& 0& 0& 0& 0& \\
1& 0& 0& 0& 0& 0& 0& 0& 0& 0& 0& \\
2& 0& 1& 0& 0& 0& 0& 0& 0& 0& 0& \\
3& 0& 1& 6& 0& 0& 0& 0& 0& 0& 0& \\
4& 0& 1& 12& 48& 0& 0& 0& 0& 0& 0& \\
5& 0& 1& 18& 156& 500& 0& 0& 0& 0& & \\
6& 0& 1& 25& 340& 2360& 6480& 0& 0& 0& & \\
7& 0& 1& 33& 636& 7060& 41400& 100842& & & & \\
8& 0& 1& 42& 1092& 17290& 162120& & & & & \\
9& 0& 1& 52& 1764& 37740& & & & & & \\
10& 0& 1& 63& 2718& & & & & & & \\
\hline
\end{tabular}
\caption{\texttt{-dc-im-l} vertex-labeled}
\label{tabL.-dc-im-l}
\end{table}

\begin{table}
\begin{tabular}{r|rrrrrrrrrrrrrrrrrr}
$E\backslash V$ & 1 & 2 & 3 & 4 & 5 & 6 & 7 &8 &9 & 10\\
\hline
0& 0& 0& 0& 0& 0& 0& 0& 0& 0& 0& \\
1& 0& 0& 0& 0& 0& 0& 0& 0& 0& 0& \\
2& 0& 0& 1& 1& 1& 1& 1& 1& 1& 1& \\
3& 0& 0& 1& 2& 3& 3& 3& 3& 3& 3& \\
4& 0& 0& 1& 4& 9& 11& 12& 12& 12& 12& \\
5& 0& 0& 1& 5& 17& 29& 37& 39& 40& & \\
6& 0& 0& 1& 7& 31& 70& 111& 132& 141& & \\
7& 0& 0& 1& 8& 48& 145& 289& & & & \\
8& 0& 0& 1& 10& 75& 289& & & & & \\
9& 0& 0& 1& 12& 111& & & & & & \\
10& 0& 0& 1& 14& & & & & & & \\
\hline
\end{tabular}
\caption{\texttt{-d-cim-l} unlabeled}
\label{tabU.-d-cim-l}
\end{table}

\begin{table}
\begin{tabular}{r|rrrrrrrrrrrrrrrrrr}
$E\backslash V$ & 1 & 2 & 3 & 4 & 5 & 6 & 7 &8 &9 & 10\\
\hline
0& 0& 0& 0& 0& 0& 0& 0& 0& 0& 0& \\
1& 0& 0& 0& 0& 0& 0& 0& 0& 0& 0& \\
2& 0& 0& 3& 6& 10& 15& 21& 28& 36& 45& \\
3& 0& 0& 3& 30& 100& 225& 441& 784& 1296& 2025& \\
4& 0& 0& 3& 54& 415& 1650& 4641& 10990& 23346& 45585& \\
5& 0& 0& 3& 78& 1030& 7215& 31521& 102676& 281016& & \\
6& 0& 0& 3& 106& 2035& 22400& 150766& 700378& 2529696& & \\
7& 0& 0& 3& 138& 3610& 56745& 557676& & & & \\
8& 0& 0& 3& 174& 5995& 127170& & & & & \\
9& 0& 0& 3& 214& 9470& & & & & & \\
10& 0& 0& 3& 258& & & & & & & \\
\hline
\end{tabular}
\caption{\texttt{-d-cim-l} vertex-labeled}
\label{tabL.-d-cim-l}
\end{table}

\begin{table}
\begin{tabular}{r|rrrrrrrrrrrrrrrrrr}
$E\backslash V$ & 1 & 2 & 3 & 4 & 5 & 6 & 7 &8 &9 & 10\\
\hline
0& 0& 0& 0& 0& 0& 0& 0& 0& 0& 0& \\
1& 0& 0& 0& 0& 0& 0& 0& 0& 0& 0& \\
2& 0& 0& 0& 0& 0& 0& 0& 0& 0& 0& \\
3& 0& 0& 0& 1& 0& 0& 0& 0& 0& 0& \\
4& 0& 0& 0& 2& 3& 1& 0& 0& 0& 0& \\
5& 0& 0& 0& 1& 7& 10& 3& 1& 0& & \\
6& 0& 0& 0& 1& 8& 28& 28& 11& 3& & \\
7& 0& 0& 0& 0& 6& 42& 91& & & & \\
8& 0& 0& 0& 0& 3& 48& & & & & \\
9& 0& 0& 0& 0& 1& & & & & & \\
10& 0& 0& 0& 0& & & & & & & \\
\hline
\end{tabular}
\caption{\texttt{-d-c-i-ml} unlabeled}
\label{tabU.-d-c-i-ml}
\end{table}

\begin{table}
\begin{tabular}{r|rrrrrrrrrrrrrrrrrr}
$E\backslash V$ & 1 & 2 & 3 & 4 & 5 & 6 & 7 &8 &9 & 10\\
\hline
0& 0& 0& 0& 0& 0& 0& 0& 0& 0& 0& \\
1& 0& 0& 0& 0& 0& 0& 0& 0& 0& 0& \\
2& 0& 0& 0& 0& 0& 0& 0& 0& 0& 0& \\
3& 0& 0& 0& 12& 0& 0& 0& 0& 0& 0& \\
4& 0& 0& 0& 18& 150& 90& 0& 0& 0& 0& \\
5& 0& 0& 0& 12& 350& 2205& 2205& 840& 0& & \\
6& 0& 0& 0& 3& 400& 6960& 37485& 49980& 34020& & \\
7& 0& 0& 0& 0& 250& 11700& 151214& & & & \\
8& 0& 0& 0& 0& 80& 12330& & & & & \\
9& 0& 0& 0& 0& 10& & & & & & \\
10& 0& 0& 0& 0& & & & & & & \\
\hline
\end{tabular}
\caption{\texttt{-d-c-i-ml} vertex-labeled}
\label{tabL.-d-c-i-ml}
\end{table}

\begin{table}
\begin{tabular}{r|rrrrrrrrrrrrrrrrrr}
$E\backslash V$ & 1 & 2 & 3 & 4 & 5 & 6 & 7 &8 &9 & 10\\
\hline
0& 0& 0& 0& 0& 0& 0& 0& 0& 0& 0& \\
1& 0& 0& 0& 0& 0& 0& 0& 0& 0& 0& \\
2& 0& 1& 0& 0& 0& 0& 0& 0& 0& 0& \\
3& 0& 1& 2& 0& 0& 0& 0& 0& 0& 0& \\
4& 0& 0& 3& 4& 0& 0& 0& 0& 0& 0& \\
5& 0& 0& 2& 10& 9& 0& 0& 0& 0& & \\
6& 0& 0& 1& 12& 30& 20& 0& 0& 0& & \\
7& 0& 0& 0& 10& 57& 93& 48& & & & \\
8& 0& 0& 0& 5& 73& 240& & & & & \\
9& 0& 0& 0& 2& 67& & & & & & \\
10& 0& 0& 0& 1& & & & & & & \\
\hline
\end{tabular}
\caption{\texttt{-dc-i-ml} unlabeled}
\label{tabU.-dc-i-ml}
\end{table}

\begin{table}
\begin{tabular}{r|rrrrrrrrrrrrrrrrrr}
$E\backslash V$ & 1 & 2 & 3 & 4 & 5 & 6 & 7 &8 &9 & 10\\
\hline
0& 0& 0& 0& 0& 0& 0& 0& 0& 0& 0& \\
1& 0& 0& 0& 0& 0& 0& 0& 0& 0& 0& \\
2& 0& 2& 0& 0& 0& 0& 0& 0& 0& 0& \\
3& 0& 1& 9& 0& 0& 0& 0& 0& 0& 0& \\
4& 0& 0& 12& 64& 0& 0& 0& 0& 0& 0& \\
5& 0& 0& 6& 156& 625& 0& 0& 0& 0& & \\
6& 0& 0& 1& 178& 2360& 7776& 0& 0& 0& & \\
7& 0& 0& 0& 116& 4495& 41400& 117649& & & & \\
8& 0& 0& 0& 45& 5495& 115020& & & & & \\
9& 0& 0& 0& 10& 4710& & & & & & \\
10& 0& 0& 0& 1& & & & & & & \\
\hline
\end{tabular}
\caption{\texttt{-dc-i-ml} vertex-labeled}
\label{tabL.-dc-i-ml}
\end{table}

\clearpage
\begin{table}
\begin{tabular}{r|rrrrrrrrrrrrrrrrrr}
$E\backslash V$ & 1 & 2 & 3 & 4 & 5 & 6 & 7 &8 &9 & 10\\
\hline
0& 0& 0& 0& 0& 0& 0& 0& 0& 0& 0& \\
1& 0& 1& 1& 1& 1& 1& 1& 1& 1& 1& \\
2& 0& 1& 3& 3& 3& 3& 3& 3& 3& 3& \\
3& 0& 0& 3& 7& 9& 9& 9& 9& 9& 9& \\
4& 0& 0& 1& 9& 20& 25& 27& 27& 27& 27& \\
5& 0& 0& 0& 6& 30& 58& 74& 79& 81& & \\
6& 0& 0& 0& 3& 32& 104& 183& 226& 243& & \\
7& 0& 0& 0& 1& 27& 149& 381& & & & \\
8& 0& 0& 0& 0& 16& 175& & & & & \\
9& 0& 0& 0& 0& 7& & & & & & \\
10& 0& 0& 0& 0& & & & & & & \\
\hline
\end{tabular}
\caption{\texttt{-d-ci-ml} unlabeled}
\label{tabU.-d-ci-ml}
\end{table}

\begin{table}
\begin{tabular}{r|rrrrrrrrrrrrrrrrrr}
$E\backslash V$ & 1 & 2 & 3 & 4 & 5 & 6 & 7 &8 &9 & 10\\
\hline
0& 0& 0& 0& 0& 0& 0& 0& 0& 0& 0& \\
1& 0& 2& 3& 4& 5& 6& 7& 8& 9& 10& \\
2& 0& 1& 12& 30& 60& 105& 168& 252& 360& 495& \\
3& 0& 0& 10& 88& 335& 875& 1946& 3864& 7050& 12045& \\
4& 0& 0& 3& 113& 1005& 4530& 14490& 38430& 90090& 192060& \\
5& 0& 0& 0& 78& 1776& 15141& 75726& 277872& 844767& & \\
6& 0& 0& 0& 28& 2035& 34523& 284991& 1521072& 6163248& & \\
7& 0& 0& 0& 4& 1570& 56745& 798897& & & & \\
8& 0& 0& 0& 0& 815& 69705& & & & & \\
9& 0& 0& 0& 0& 275& & & & & & \\
10& 0& 0& 0& 0& & & & & & & \\
\hline
\end{tabular}
\caption{\texttt{-d-ci-ml} vertex-labeled}
\label{tabL.-d-ci-ml}
\end{table}

\begin{table}
\begin{tabular}{r|rrrrrrrrrrrrrrrrrr}
$E\backslash V$ & 1 & 2 & 3 & 4 & 5 & 6 & 7 &8 &9 & 10\\
\hline
0& 0& 0& 0& 0& 0& 0& 0& 0& 0& 0& \\
1& 0& 0& 0& 0& 0& 0& 0& 0& 0& 0& \\
2& 0& 0& 0& 0& 0& 0& 0& 0& 0& 0& \\
3& 0& 0& 0& 0& 0& 0& 0& 0& 0& 0& \\
4& 0& 0& 0& 3& 0& 0& 0& 0& 0& 0& \\
5& 0& 0& 0& 11& 10& 3& 0& 0& 0& & \\
6& 0& 0& 0& 27& 51& 44& 11& 3& 0& & \\
7& 0& 0& 0& 51& 157& 236& 153& & & & \\
8& 0& 0& 0& 93& 386& 850& & & & & \\
9& 0& 0& 0& 150& 838& & & & & & \\
10& 0& 0& 0& 241& & & & & & & \\
\hline
\end{tabular}
\caption{\texttt{-d-c-iml} unlabeled}
\label{tabU.-d-c-iml}
\end{table}

\begin{table}
\begin{tabular}{r|rrrrrrrrrrrrrrrrrr}
$E\backslash V$ & 1 & 2 & 3 & 4 & 5 & 6 & 7 &8 &9 & 10\\
\hline
0& 0& 0& 0& 0& 0& 0& 0& 0& 0& 0& \\
1& 0& 0& 0& 0& 0& 0& 0& 0& 0& 0& \\
2& 0& 0& 0& 0& 0& 0& 0& 0& 0& 0& \\
3& 0& 0& 0& 0& 0& 0& 0& 0& 0& 0& \\
4& 0& 0& 0& 36& 0& 0& 0& 0& 0& 0& \\
5& 0& 0& 0& 144& 600& 360& 0& 0& 0& & \\
6& 0& 0& 0& 360& 3250& 11925& 11025& 4200& 0& & \\
7& 0& 0& 0& 738& 10650& 76635& 257985& & & & \\
8& 0& 0& 0& 1365& 27650& 308385& & & & & \\
9& 0& 0& 0& 2352& 62940& & & & & & \\
10& 0& 0& 0& 3834& & & & & & & \\
\hline
\end{tabular}
\caption{\texttt{-d-c-iml} vertex-labeled}
\label{tabL.-d-c-iml}
\end{table}

\begin{table}
\begin{tabular}{r|rrrrrrrrrrrrrrrrrr}
$E\backslash V$ & 1 & 2 & 3 & 4 & 5 & 6 & 7 &8 &9 & 10\\
\hline
0& 0& 0& 0& 0& 0& 0& 0& 0& 0& 0& \\
1& 0& 0& 0& 0& 0& 0& 0& 0& 0& 0& \\
2& 0& 0& 0& 0& 0& 0& 0& 0& 0& 0& \\
3& 0& 2& 0& 0& 0& 0& 0& 0& 0& 0& \\
4& 0& 5& 5& 0& 0& 0& 0& 0& 0& 0& \\
5& 0& 8& 19& 13& 0& 0& 0& 0& 0& & \\
6& 0& 11& 45& 70& 35& 0& 0& 0& 0& & \\
7& 0& 15& 87& 227& 245& 95& 0& & & & \\
8& 0& 19& 153& 579& 1029& 840& & & & & \\
9& 0& 24& 252& 1302& 3346& & & & & & \\
10& 0& 29& 394& 2681& & & & & & & \\
\hline
\end{tabular}
\caption{\texttt{-dc-iml} unlabeled}
\label{tabU.-dc-iml}
\end{table}

\begin{table}
\begin{tabular}{r|rrrrrrrrrrrrrrrrrr}
$E\backslash V$ & 1 & 2 & 3 & 4 & 5 & 6 & 7 &8 &9 & 10\\
\hline
0& 0& 0& 0& 0& 0& 0& 0& 0& 0& 0& \\
1& 0& 0& 0& 0& 0& 0& 0& 0& 0& 0& \\
2& 0& 0& 0& 0& 0& 0& 0& 0& 0& 0& \\
3& 0& 4& 0& 0& 0& 0& 0& 0& 0& 0& \\
4& 0& 9& 27& 0& 0& 0& 0& 0& 0& 0& \\
5& 0& 14& 102& 256& 0& 0& 0& 0& 0& & \\
6& 0& 20& 240& 1420& 3125& 0& 0& 0& 0& & \\
7& 0& 27& 471& 4688& 23535& 46656& 0& & & & \\
8& 0& 35& 840& 12250& 102900& 453096& & & & & \\
9& 0& 44& 1400& 28080& 345730& & & & & & \\
10& 0& 54& 2214& 58914& & & & & & & \\
\hline
\end{tabular}
\caption{\texttt{-dc-iml} vertex-labeled}
\label{tabL.-dc-iml}
\end{table}

\begin{table}
\begin{tabular}{r|rrrrrrrrrrrrrrrrrr}
$E\backslash V$ & 1 & 2 & 3 & 4 & 5 & 6 & 7 &8 &9 & 10\\
\hline
0& 0& 0& 0& 0& 0& 0& 0& 0& 0& 0& \\
1& 0& 0& 0& 0& 0& 0& 0& 0& 0& 0& \\
2& 0& 1& 1& 1& 1& 1& 1& 1& 1& 1& \\
3& 0& 2& 6& 6& 6& 6& 6& 6& 6& 6& \\
4& 0& 3& 15& 24& 29& 29& 29& 29& 29& 29& \\
5& 0& 3& 26& 66& 107& 122& 127& 127& 127& & \\
6& 0& 4& 40& 142& 318& 454& 520& 536& 541& & \\
7& 0& 4& 57& 269& 800& 1464& 1967& & & & \\
8& 0& 5& 79& 474& 1813& 4224& & & & & \\
9& 0& 5& 106& 793& 3810& & & & & & \\
10& 0& 6& 138& 1273& & & & & & & \\
\hline
\end{tabular}
\caption{\texttt{-d-ciml} unlabeled}
\label{tabU.-d-ciml}
\end{table}

\begin{table}
\begin{tabular}{r|rrrrrrrrrrrrrrrrrr}
$E\backslash V$ & 1 & 2 & 3 & 4 & 5 & 6 & 7 &8 &9 & 10\\
\hline
0& 0& 0& 0& 0& 0& 0& 0& 0& 0& 0& \\
1& 0& 0& 0& 0& 0& 0& 0& 0& 0& 0& \\
2& 0& 2& 3& 4& 5& 6& 7& 8& 9& 10& \\
3& 0& 4& 27& 64& 125& 216& 343& 512& 729& 1000& \\
4& 0& 5& 69& 358& 1190& 2946& 6349& 12356& 22239& 37630& \\
5& 0& 6& 123& 1104& 6275& 23796& 70315& 177920& 404109& & \\
6& 0& 7& 193& 2554& 22585& 130286& 543837& 1813568& 5197044& & \\
7& 0& 8& 285& 5102& 64340& 538614& 3165841& & & & \\
8& 0& 9& 402& 9363& 158520& 1829799& & & & & \\
9& 0& 10& 547& 16176& 354905& & & & & & \\
10& 0& 11& 723& 26626& & & & & & & \\
\hline
\end{tabular}
\caption{\texttt{-d-ciml} vertex-labeled}
\label{tabL.-d-ciml}
\end{table}

\clearpage
\section{Accumulated Marginal statistics}
\label{sec.marg}

Adding the contents of one or more of the
previous arrays defines the union of their graph sets,
and regards some of the properties as irrelevant in these tables.
If we look on the flags as defining a hypertable along five or six axes,
these sums are the marginal sums; they create the Tables in Section \ref{sec.marg}.

The properties that are not taken into account while
counting the graphs are either replaced by the filler \texttt{.*}
or not tagged at all, using regular expressions
of the usual programming languages as the tags.
\begin{exa}
\texttt{-d.*-m-l} flags graphs that are undirected, have any type of isolated vertices
or connectivity, but have no multiedges or loops.
\end{exa}
\begin{exa}
\texttt{-dc} flags graphs that are undirected and connected, but have
any type of isolated vertices, multiedges or loops.
\end{exa}

\subsection{Undirected Graphs}
Tables \ref{tab.046742}--\ref{tabL.-d-l}
summarize statistics of undirected graphs.

\begin{table}
\begin{tabular}{r|rrrrrrrrrrrrrrr}
$E\backslash V$ & 1 & 2 & 3 & 4 & 5 & 6 & 7 &8 &9\\
\hline
0 &1\\
1 &0 &1\\
2 &0 &0 &1\\
3 &0 &0 &1 &2\\
4 &0 &0 &0 &2 &3\\
5 &0 &0 &0 &1 &5 &6\\
6 &0 &0 &0 &1 &5 &13 &11\\
7 &0 &0 &0 &0 &4 &19 &33 &23\\
8 &0 &0 &0 &0 &2 &22 &67 &89 &47\\
9 &0 &0 &0 &0 &1 &20 &107 &236 &240 &106\\
10 &0 &0 &0 &0 &1 &14 &132 &486 &797 &657 &235\\
\hline
\end{tabular}

\caption{\texttt{( -dc.*-m-l)} unlabeled.
The number of connected undirected graphs without multiedges
or loops \cite[A054923,A054924,A046742]{EIS}\cite[Vol. 1, Sec. 7, Table 1]{SteinbachVol4}.
With the exception of 
the value of 1 at $E=0$, $V=1$ this is the same as
Table \ref{tabU.-dc-i-m-l}.
}
\label{tab.046742}
\end{table}

\begin{table}
\begin{tabular}{r|rrrrrrrrrrrrrrrrrr}
$E\backslash V$ & 1 & 2 & 3 & 4 & 5 & 6 & 7 &8 &9 & 10\\
\hline
0& 1& 0& 0& 0& 0& 0& 0& 0& 0& 0& & & \\
1& 0& 1& 0& 0& 0& 0& 0& 0& 0& 0& & & \\
2& 0& 0& 3& 0& 0& 0& 0& 0& 0& 0& & & \\
3& 0& 0& 1& 16& 0& 0& 0& 0& 0& 0& & & \\
4& 0& 0& 0& 15& 125& 0& 0& 0& 0& 0& & & \\
5& 0& 0& 0& 6& 222& 1296& 0& 0& 0& & & & \\
6& 0& 0& 0& 1& 205& 3660& 16807& 0& & & & & \\
7& 0& 0& 0& 0& 120& 5700& & & & & & & \\
8& 0& 0& 0& 0& 45& & & & & & & & \\
9& 0& 0& 0& 0& & & & & & & & & \\
\hline
\end{tabular}

\caption{\texttt{( -dc.*-m-l)} vertex-labeled \cite[A062734]{EIS}.}
\end{table}

\begin{table}
\begin{tabular}{r|rrrrrrrrrrrrrrrrrr}
$E\backslash V$ & 1 & 2 & 3 & 4 & 5 & 6 & 7 &8 &9 & 10\\
\hline
0& 1& 1& 1& 1& 1& 1& 1& 1& 1& 1& \\
1& 1& 2& 2& 2& 2& 2& 2& 2& 2& 2& \\
2& 1& 4& 6& 7& 7& 7& 7& 7& 7& 7& \\
3& 1& 6& 14& 20& 22& 23& 23& 23& 23& 23& \\
4& 1& 9& 28& 53& 69& 76& 78& 79& 79& 79& \\
5& 1& 12& 52& 125& 198& 245& 264& 271& 273& & \\
6& 1& 16& 93& 287& 550& 782& 915& 973& 993& & \\
7& 1& 20& 152& 606& 1441& 2392& 3111& & & & \\
8& 1& 25& 242& 1226& 3611& 7118& & & & & \\
9& 1& 30& 370& 2358& 8608& & & & & & \\
10& 1& 36& 546& 4356& & & & & & & \\
\hline
\end{tabular}

\caption{
\texttt{( -d)} unlabeled.
The number of undirected graphs 
allowing loops and multiedges
\cite[A290428]{EIS}.
Sum of Table \ref{tabU.-dc} and Table \ref{tabU.-d-c}.}
\label{tab.290428}
\end{table}

\begin{table}
\begin{tabular}{r|rrrrrrrrrrrrrrrrrr}
$E\backslash V$ & 1 & 2 & 3 & 4 & 5 & 6 & 7 &8 &9 \\
\hline
0& 1& 1& 1& 1& 1& 1& 1& 1& 1& \\
1& 1& 3& 6& 10& 15& 21& 28& 36& 45\\
2& 1& 6& 21& 55& 120& 231& 406& 666& 1035\\
3& 1& 10& 56& 220& 680& 1771& 4060& 8436& 16215\\
4& 1& 15& 126& 715& 3060& 10626& 31465& 82251& 194580\\
5& 1& 21& 252& 2002& 11628& 53130& 201376& 658008& 1906884& & \\
6& 1& 28& 462& 5005& 38760& 230230& 1107568& 4496388& 15890700& & \\
7& 1& 36& 792& 11440& 116280& 888030& 5379616& & & & \\
8& 1& 45& 1287& 24310& 319770& 3108105& & & & & \\
9& 1& 55& 2002& 48620& 817190& & & & & & \\
10& 1& 66& 3003& 92378& & & & & & & \\
\hline
\end{tabular}

\caption{\texttt{( -d)} vertex-labeled  \cite[A098568]{EIS}.
Sum of Table \ref{tabL.-dc} and Table \ref{tabL.-d-c}.}
\end{table}

\begin{table}
\begin{tabular}{r|rrrrrrrrrrrrrrrrrr}
$E\backslash V$ & 1 & 2 & 3 & 4 & 5 & 6 & 7 &8 &9 & 10\\
\hline
0 &1\\
1 &0 &1\\
2 &0 &1 &1\\
3 &0 &1 &2 &2\\
4 &0 &1 &3 &5 &3\\
5 &0 &1 &4 &11 &11 &6\\
6 &0 &1 &6 &22 &34 &29 &11\\
7 &0 &1 &7 &37 &85 &110 &70 &23\\
8 &0 &1 &9 &61 &193 &348 &339 &185 &47\\
9 &0 &1 &11 &95 &396 &969 &1318 &1067 &479 &106\\
10 &0 &1 &13 &141 &771 &2445 &4457 &4940 &3294 &1729 &235\\
\hline
\end{tabular}

\caption{ \texttt{-dc.*-l} unlabeled.
Undirected loopless connected multigraphs with $E$ edges and $V$ vertices \cite[A191646]{EIS}.
With the exception of the 1 at $E=0$, $V=1$,
the sum of Table \ref{tabU.-dc-i-m-l} and Table \ref{tabU.-dc-im-l}.}
\label{tab.191646}
\end{table}

\begin{table}
\begin{tabular}{r|rrrrrrrrrrrrrrrrrr}
$E\backslash V$ & 1 & 2 & 3 & 4 & 5 & 6 & 7 &8 &9 & 10\\
\hline
0& 1& 0& 0& 0& 0& 0& 0& 0& 0& 0& \\
1& 0& 1& 0& 0& 0& 0& 0& 0& 0& 0& \\
2& 0& 1& 3& 0& 0& 0& 0& 0& 0& 0& \\
3& 0& 1& 7& 16& 0& 0& 0& 0& 0& 0& \\
4& 0& 1& 12& 63& 125& 0& 0& 0& 0& 0& \\
5& 0& 1& 18& 162& 722& 1296& 0& 0& 0& & \\
6& 0& 1& 25& 341& 2565& 10140& 16807& 0& 0& & \\
7& 0& 1& 33& 636& 7180& 47100& 169137& & & & \\
8& 0& 1& 42& 1092& 17335& 168285& & & & & \\
9& 0& 1& 52& 1764& 37750& & & & & & \\
10& 0& 1& 63& 2718& & & & & & & \\
\hline
\end{tabular}

\caption{\texttt{( -dc.*-l)} vertex-labeled \cite[A290776]{EIS}.
}
\end{table}

\begin{table}
\begin{tabular}{r|rrrrrrrrrrrrrrrrrr}
$E\backslash V$ & 1 & 2 & 3 & 4 & 5 & 6 & 7 &8 &9 & 10\\
\hline
0& 1& 0& 0& 0& 0& 0& 0& 0& 0& 0& \\
1& 1& 1& 0& 0& 0& 0& 0& 0& 0& 0& \\
2& 1& 2& 1& 0& 0& 0& 0& 0& 0& 0& \\
3& 1& 4& 4& 2& 0& 0& 0& 0& 0& 0& \\
4& 1& 6& 11& 9& 3& 0& 0& 0& 0& 0& \\
5& 1& 9& 25& 34& 20& 6& 0& 0& 0& & \\
6& 1& 12& 52& 104& 99& 49& 11& 0& 0& & \\
7& 1& 16& 94& 274& 387& 298& 118& & & & \\
8& 1& 20& 162& 645& 1295& 1428& & & & & \\
9& 1& 25& 263& 1399& 3809& & & & & & \\
10& 1& 30& 407& 2823& & & & & & & \\
\hline
\end{tabular}

\caption{\texttt{( -dc)} unlabeled.
Row sums in \cite[A007719]{EIS}
}
\label{tabU.-dc}
\end{table}

\begin{table}
\begin{tabular}{r|rrrrrrrrrrrrrrrrrr}
$E\backslash V$ & 1 & 2 & 3 & 4 & 5 & 6 & 7 &8 &9 & 10\\
\hline
0& 1& 0& 0& 0& 0& 0& 0& 0& 0& 0& \\
1& 1& 1& 0& 0& 0& 0& 0& 0& 0& 0& \\
2& 1& 3& 3& 0& 0& 0& 0& 0& 0& 0& \\
3& 1& 6& 16& 16& 0& 0& 0& 0& 0& 0& \\
4& 1& 10& 51& 127& 125& 0& 0& 0& 0& 0& \\
5& 1& 15& 126& 574& 1347& 1296& 0& 0& 0& & \\
6& 1& 21& 266& 1939& 8050& 17916& 16807& 0& 0& & \\
7& 1& 28& 504& 5440& 35210& 135156& 286786& & & & \\
8& 1& 36& 882& 13387& 125730& 736401& & & & & \\
9& 1& 45& 1452& 29854& 388190& & & & & & \\
10& 1& 55& 2277& 61633& & & & & & & \\
\hline
\end{tabular}

\caption{\texttt{( -dc)} vertex-labeled.}
\label{tabL.-dc}
\end{table}

\begin{table}
\begin{tabular}{r|rrrrrrrrrrrrrrr}
$E\backslash V$ & 1 & 2 & 3 & 4 & 5 & 6 & 7 &8 &9 &10 &11 & 12\\
\hline
0 &1 &1 &1 &1 & 1 & 1 & 1 & 1 & 1 & 1 & 1\\
1 &0 &1 &1 &1 &1 & 1 & 1 & 1 & 1 & 1 & 1\\
2 &0 &0 &1 &2 &2 &2& 2& 2& 2& 2& 2\\
3 &0 &0 &1 &3 &4 &5 &5 & 5 & 5 & 5 & 5 \\
4 &0 &0 &0 &2 &6 &9 &10 &11 & 11 & 11 & 11 \\
5 &0 &0 &0 &1 &6 &15 &21 &24 &25 & 26 & 26 \\
6 &0 &0 &0 &1 &6 &21 &41 &56 &63 &66 &67 \\
7 &0 &0 &0 &0 &4 &24 &65 &115 &148 &165 &172\\
8 &0 &0 &0 &0 &2 &24 &97 &221 &345 &428 &467 \\
9 &0 &0 &0 &0 &1 &21 &131 &402 &771 &1103 &1305 &1405\\
10 &0 &0 &0 &0 &1 &15 &148 &663 &1637 &2769 &3664 &4191\\
\end{tabular}
\caption{\texttt{( -d.*-m-l)} unlabeled.
Simple graphs with $E$ edges and $V$ vertices \cite[A008406]{EIS}\cite[vol. 4, Tables 2.2--2.2g]{SteinbachVol4}.
Sum of tables \ref{tab.046742}, \ref{tab.sk2}--\ref{tab.sk5} and contributions
by more than 5 components.
}
\label{tab.008406}
\end{table}

\begin{table}
\begin{tabular}{r|rrrrrrrrrrrrrrr}
$E\backslash V$ & 1 & 2 & 3 & 4 & 5 & 6 & 7 &8 &9 & 10\\
\hline
0& 1& 1& 1& 1& 1& 1& 1& 1& 1& 1& & & & & \\
1& 0& 1& 3& 6& 10& 15& 21& 28& 36& 45& & & & & \\
2& 0& 0& 3& 15& 45& 105& 210& 378& 630& 990& & & & & \\
3& 0& 0& 1& 20& 120& 455& 1330& 3276& 7140& 14190& & & & & \\
4& 0& 0& 0& 15& 210& 1365& 5985& 20475& 58905& & & & & & \\
5& 0& 0& 0& 6& 252& 3003& 20349& 98280& 376992& & & & & & \\
6& 0& 0& 0& 1& 210& 5005& 54264& 376740& & & & & & & \\
7& 0& 0& 0& 0& 120& 6435& & & & & & & & & \\
8& 0& 0& 0& 0& 45& & & & & & & & & & \\
9& 0& 0& 0& & & & & & & & & & & & \\
\hline
\end{tabular}

\caption{\texttt{( -d.*-m-l)} vertex-labeled  \cite[A084546]{EIS}.
}
\end{table}

\begin{table}
\begin{tabular}{r|rrrrrrrrrrrrrrrrrr}
$E\backslash V$ & 1 & 2 & 3 & 4 & 5 & 6 & 7 &8 &9 & 10\\
\hline
0& 0& 1& 1& 1& 1& 1& 1& 1& 1& 1& \\
1& 0& 1& 2& 2& 2& 2& 2& 2& 2& 2& \\
2& 0& 2& 5& 7& 7& 7& 7& 7& 7& 7& \\
3& 0& 2& 10& 18& 22& 23& 23& 23& 23& 23& \\
4& 0& 3& 17& 44& 66& 76& 78& 79& 79& 79& \\
5& 0& 3& 27& 91& 178& 239& 264& 271& 273& & \\
6& 0& 4& 41& 183& 451& 733& 904& 973& 993& & \\
7& 0& 4& 58& 332& 1054& 2094& 2993& & & & \\
8& 0& 5& 80& 581& 2316& 5690& & & & & \\
9& 0& 5& 107& 959& 4799& & & & & & \\
10& 0& 6& 139& 1533& & & & & & & \\
\hline
\end{tabular}

\caption{\texttt{( -d-c)} unlabeled.
See \cite[A007717]{EIS} for the limit $V\to\infty$.
}
\label{tabU.-d-c}
\end{table}

\begin{table}
\begin{tabular}{r|rrrrrrrrrrrrrrrrrr}
$E\backslash V$ & 1 & 2 & 3 & 4 & 5 & 6 & 7 &8 &9 \\
\hline
0& 0& 1& 1& 1& 1& 1& 1& 1& 1& \\
1& 0& 2& 6& 10& 15& 21& 28& 36& 45\\
2& 0& 3& 18& 55& 120& 231& 406& 666& 1035\\
3& 0& 4& 40& 204& 680& 1771& 4060& 8436& 16215\\
4& 0& 5& 75& 588& 2935& 10626& 31465& 82251& 194580\\
5& 0& 6& 126& 1428& 10281& 51834& 201376& 658008& 1906884& & \\
6& 0& 7& 196& 3066& 30710& 212314& 1090761& 4496388& 15890700& & \\
7& 0& 8& 288& 6000& 81070& 752874& 5092830& & & & \\
8& 0& 9& 405& 10923& 194040& 2371704& & & & & \\
9& 0& 10& 550& 18766& 429000& & & & & & \\
10& 0& 11& 726& 30745& & & & & & & \\
\hline
\end{tabular}

\caption{\texttt{( -d-c)} vertex-labeled.}
\label{tabL.-d-c}
\end{table}

\begin{table}
\begin{tabular}{r|rrrrrrrrrrrrrrrrrr}
$E\backslash V$ & 1 & 2 & 3 & 4 & 5 & 6 & 7 &8 &9 & 10\\
\hline
0& 1& 1& 1& 1& 1& 1& 1& 1& 1& 1& \\
1& 0& 1& 1& 1& 1& 1& 1& 1& 1& 1& \\
2& 0& 1& 2& 3& 3& 3& 3& 3& 3& 3& \\
3& 0& 1& 3& 6& 7& 8& 8& 8& 8& 8& \\
4& 0& 1& 4& 11& 17& 21& 22& 23& 23& 23& \\
5& 0& 1& 5& 18& 35& 52& 60& 64& 65& & \\
6& 0& 1& 7& 32& 76& 132& 173& 197& 206& & \\
7& 0& 1& 8& 48& 149& 313& 471& & & & \\
8& 0& 1& 10& 75& 291& 741& & & & & \\
9& 0& 1& 12& 111& 539& & & & & & \\
10& 0& 1& 14& 160& & & & & & & \\
\hline
\end{tabular}

\caption{\texttt{( -d.*-l)} unlabeled.
Undirected loopless multigraphs with $E$ edges and $V$ vertices \cite[A192517]{EIS}.
Sum of tables \ref{tab.191646}, \ref{tab.mk2}--\ref{tab.mk5} and so forth.
}
\label{tab.192517}
\end{table}

\begin{table}
\begin{tabular}{r|rrrrrrrrrrrrrrrrrr}
$E\backslash V$ & 1 & 2 & 3 & 4 & 5 & 6 & 7 &8 &9 & 10\\
\hline
0& 1& 1& 1& 1& 1& 1& 1& 1& 1& 1& \\
1& 0& 1& 3& 6& 10& 15& 21& 28& 36& 45& \\
2& 0& 1& 6& 21& 55& 120& 231& 406& 666& 1035& \\
3& 0& 1& 10& 56& 220& 680& 1771& 4060& 8436& 16215& \\
4& 0& 1& 15& 126& 715& 3060& 10626& 31465& 82251& 194580& \\
5& 0& 1& 21& 252& 2002& 11628& 53130& 201376& 658008& & \\
6& 0& 1& 28& 462& 5005& 38760& 230230& 1107568& 4496388& & \\
7& 0& 1& 36& 792& 11440& 116280& 888030& & & & \\
8& 0& 1& 45& 1287& 24310& 319770& & & & & \\
9& 0& 1& 55& 2002& 48620& & & & & & \\
10& 0& 1& 66& 3003& & & & & & & \\
\hline
\end{tabular}

\caption{\texttt{( -d.*-l)} vertex-labeled.
}
\label{tabL.-d-l}
\end{table}

\clearpage
\subsection{Directed Graphs}
Tables \ref{tab.139621}--\ref{tab.fin}
summarize statistics of oriented/directed graphs.

\begin{table}
\begin{tabular}{r|rrrrrrrrrrrrrrrrrr}
$E\backslash V$ & 1 & 2 & 3 & 4 & 5 & 6 & 7 &8 &9 & 10\\
\hline
0& 0& 0& 0& 0& 0& 0& 0& 0& 0& 0& \\
1& 0& 1& 0& 0& 0& 0& 0& 0& 0& 0& \\
2& 0& 4& 3& 0& 0& 0& 0& 0& 0& 0& \\
3& 0& 8& 15& 8& 0& 0& 0& 0& 0& 0& \\
4& 0& 16& 57& 66& 27& 0& 0& 0& 0& 0& \\
5& 0& 25& 163& 353& 295& 91& 0& 0& 0& & \\
6& 0& 40& 419& 1504& 2203& 1407& 350& 0& 0& & \\
7& 0& 56& 932& 5302& 12382& 13372& 6790& & & & \\
8& 0& 80& 1940& 16549& 58237& 96456& & & & & \\
9& 0& 105& 3743& 46566& 237904& & & & & & \\
10& 0& 140& 6867& 121111& & & & & & & \\
\hline
\end{tabular}

\caption{
\texttt{( d.*Cc-i)} unlabeled.
The number of connected directed multigraphs with loops
and no isolated vertex, with $E$ arcs and $V$ vertices
\cite[A139621]{EIS}.
}
\label{tab.139621}
\end{table}

\begin{table}
\begin{tabular}{r|rrrrrrrrrrrrrrrrrr}
$E\backslash V$ & 1 & 2 & 3 & 4 & 5 & 6 & 7 &8 &9 & 10\\
\hline
0& 0& 0& 0& 0& 0& 0& 0& 0& 0& 0& \\
1& 0& 2& 0& 0& 0& 0& 0& 0& 0& 0& \\
2& 0& 7& 12& 0& 0& 0& 0& 0& 0& 0& \\
3& 0& 16& 80& 128& 0& 0& 0& 0& 0& 0& \\
4& 0& 30& 315& 1328& 2000& 0& 0& 0& 0& 0& \\
5& 0& 50& 951& 7808& 29104& 41472& 0& 0& 0& & \\
6& 0& 77& 2429& 34136& 234920& 794112& 1075648& 0& 0& & \\
7& 0& 112& 5517& 123272& 1386880& 8328192& 25952128& & & & \\
8& 0& 156& 11475& 388223& 6674205& 63248832& & & & & \\
9& 0& 210& 22275& 1101408& 27706645& & & & & & \\
10& 0& 275& 40887& 2875224& & & & & & & \\
\hline
\end{tabular}

\caption{\texttt{( d.*Cc-i)} vertex-labeled}
\end{table}

\begin{table}
\begin{tabular}{r|rrrrrrrrrrrrrrrrrr}
$E\backslash V$ & 1 & 2 & 3 & 4 & 5 & 6 & 7 &8 &9 & 10\\
\hline
0& 1& 0& 0& 0& 0& 0& 0& 0& 0& 0& \\
1& 1& 0& 0& 0& 0& 0& 0& 0& 0& 0& \\
2& 1& 1& 0& 0& 0& 0& 0& 0& 0& 0& \\
3& 1& 2& 1& 0& 0& 0& 0& 0& 0& 0& \\
4& 1& 6& 4& 1& 0& 0& 0& 0& 0& 0& \\
5& 1& 10& 19& 6& 1& 0& 0& 0& 0& & \\
6& 1& 19& 73& 59& 9& 1& 0& 0& 0& & \\
7& 1& 28& 208& 350& 138& 12& 1& & & & \\
8& 1& 44& 534& 1670& 1361& 301& & & & & \\
9& 1& 60& 1215& 6476& 9724& & & & & & \\
10& 1& 85& 2542& 21898& & & & & & & \\
\hline
\end{tabular}

\caption{
\texttt{( dC)} unlabeled.
The number of strongly connected directed multigraphs with loops
and no vertex of degree zero, with n arcs and k vertices
\cite[A139622]{EIS}.
}
\label{tab.139622}
\end{table}

\begin{table}
\begin{tabular}{r|rrrrrrrrrrrrrrrrrr}
$E\backslash V$ & 1 & 2 & 3 & 4 & 5 & 6 & 7 &8 &9 & 10\\
\hline
0& 1& 0& 0& 0& 0& 0& 0& 0& 0& 0& \\
1& 1& 0& 0& 0& 0& 0& 0& 0& 0& 0& \\
2& 1& 1& 0& 0& 0& 0& 0& 0& 0& 0& \\
3& 1& 4& 2& 0& 0& 0& 0& 0& 0& 0& \\
4& 1& 10& 21& 6& 0& 0& 0& 0& 0& 0& \\
5& 1& 20& 111& 132& 24& 0& 0& 0& 0& & \\
6& 1& 35& 413& 1288& 960& 120& 0& 0& 0& & \\
7& 1& 56& 1233& 8152& 15680& 7920& 720& & & & \\
8& 1& 84& 3159& 39049& 156955& 201450& & & & & \\
9& 1& 120& 7227& 153540& 1140055& & & & & & \\
10& 1& 165& 15147& 520404& & & & & & & \\
\hline
\end{tabular}

\caption{
\texttt{( dC)} vertex-labeled.
}
\end{table}

\begin{table}
\begin{tabular}{r|rrrrrrrrrrrrrrrrrr}
$E\backslash V$ & 1 & 2 & 3 & 4 & 5 & 6 & 7 &8 &9 & 10\\
\hline
0& 1& 0& 0& 0& 0& 0& 0& 0& 0& 0& \\
1& 0& 1& 0& 0& 0& 0& 0& 0& 0& 0& \\
2& 0& 1& 3& 0& 0& 0& 0& 0& 0& 0& \\
3& 0& 0& 4& 8& 0& 0& 0& 0& 0& 0& \\
4& 0& 0& 4& 22& 27& 0& 0& 0& 0& 0& \\
5& 0& 0& 1& 37& 108& 91& 0& 0& 0& & \\
6& 0& 0& 1& 47& 326& 582& 350& 0& 0& & \\
7& 0& 0& 0& 38& 667& 2432& 3024& & & & \\
8& 0& 0& 0& 27& 1127& 7694& & & & & \\
9& 0& 0& 0& 13& 1477& & & & & & \\
10& 0& 0& 0& 5& & & & & & & \\
\hline
\end{tabular}

\caption{
\texttt{( d.*Cc.*-m-l)} unlabeled.
The number of weakly connected directed graphs without multiedges
or loops \cite[A054733,A283753]{EIS}.
The undirected variants are 
in Table \ref{tab.046742}.
}
\label{tab.054733}
\end{table}

\begin{table}
\begin{tabular}{r|rrrrrrrrrrrrrrrrrr}
$E\backslash V$ & 1 & 2 & 3 & 4 & 5 & 6 & 7 &8 &9 & 10\\
\hline
0& 1& 0& 0& 0& 0& 0& 0& 0& 0& 0& \\
1& 0& 2& 0& 0& 0& 0& 0& 0& 0& 0& \\
2& 0& 1& 12& 0& 0& 0& 0& 0& 0& 0& \\
3& 0& 0& 20& 128& 0& 0& 0& 0& 0& 0& \\
4& 0& 0& 15& 432& 2000& 0& 0& 0& 0& 0& \\
5& 0& 0& 6& 768& 11104& 41472& 0& 0& 0& & \\
6& 0& 0& 1& 920& 33880& 337920& 1075648& 0& 0& & \\
7& 0& 0& 0& 792& 73480& 1536000& 11968704& & & & \\
8& 0& 0& 0& 495& 123485& 5062080& & & & & \\
9& 0& 0& 0& 220& 166860& & & & & & \\
10& 0& 0& 0& 66& & & & & & & \\
\hline
\end{tabular}

\caption{\texttt{( d.*Cc.*-m-l)} vertex-labeled.
\cite[A062735]{EIS}
}
\end{table}

\begin{table}
\begin{tabular}{r|rrrrrrrrrrrrrrrrrr}
$E\backslash V$ & 1 & 2 & 3 & 4 & 5 & 6 & 7 &8 &9 & 10\\
\hline
0& 1& 0& 0& 0& 0& 0& 0& 0& 0& 0& \\
1& 0& 0& 0& 0& 0& 0& 0& 0& 0& 0& \\
2& 0& 1& 0& 0& 0& 0& 0& 0& 0& 0& \\
3& 0& 0& 1& 0& 0& 0& 0& 0& 0& 0& \\
4& 0& 0& 2& 1& 0& 0& 0& 0& 0& 0& \\
5& 0& 0& 1& 4& 1& 0& 0& 0& 0& & \\
6& 0& 0& 1& 16& 7& 1& 0& 0& 0& & \\
7& 0& 0& 0& 22& 58& 10& 1& & & & \\
8& 0& 0& 0& 22& 240& 165& & & & & \\
9& 0& 0& 0& 11& 565& & & & & & \\
10& 0& 0& 0& 5& & & & & & & \\
\hline
\end{tabular}

\caption{
\texttt{( dCc.*-m-l)} unlabeled.
The number of strongly connected directed graphs without loops or multiedges.
Strongly connected variant of Table \ref{tab.054733}.
With the exception of the 1 at $E=0$, $V=1$ this is the same
as Table \ref{tabU.dCc-i-m-l}.
}
\end{table}

\begin{table}
\begin{tabular}{r|rrrrrrrrrrrrrrrrrr}
$E\backslash V$ & 1 & 2 & 3 & 4 & 5 & 6 & 7 &8 &9 & 10\\
\hline
0& 1& 0& 0& 0& 0& 0& 0& 0& 0& 0& \\
1& 0& 0& 0& 0& 0& 0& 0& 0& 0& 0& \\
2& 0& 1& 0& 0& 0& 0& 0& 0& 0& 0& \\
3& 0& 0& 2& 0& 0& 0& 0& 0& 0& 0& \\
4& 0& 0& 9& 6& 0& 0& 0& 0& 0& 0& \\
5& 0& 0& 6& 84& 24& 0& 0& 0& 0& & \\
6& 0& 0& 1& 316& 720& 120& 0& 0& 0& & \\
7& 0& 0& 0& 492& 6440& 6480& 720& & & & \\
8& 0& 0& 0& 417& 26875& 107850& & & & & \\
9& 0& 0& 0& 212& 65280& & & & & & \\
10& 0& 0& 0& 66& & & & & & & \\
\hline
\end{tabular}

\caption{\texttt{( dCc.*-m-l)} vertex-labeled.
With the exception of the 1 at $E=0$, $V=1$ this is the same
as Table \ref{tabL.dCc-i-m-l}.
}
\end{table}

\begin{table}
\begin{tabular}{r|rrrrrrrrrrrrrrrrrr}
$E\backslash V$ & 1 & 2 & 3 & 4 & 5 & 6 & 7 &8 &9 & 10\\
\hline
0& 1& 1& 1& 1& 1& 1& 1& 1& 1& 1& \\
1& 1& 2& 2& 2& 2& 2& 2& 2& 2& 2& \\
2& 1& 6& 10& 11& 11& 11& 11& 11& 11& 11& \\
3& 1& 10& 31& 47& 51& 52& 52& 52& 52& 52& \\
4& 1& 19& 90& 198& 269& 291& 295& 296& 296& 296& \\
5& 1& 28& 222& 713& 1270& 1596& 1697& 1719& 1723& & \\
6& 1& 44& 520& 2423& 5776& 8838& 10425& 10922& 11033& & \\
7& 1& 60& 1090& 7388& 24032& 46384& 63419& & & & \\
8& 1& 85& 2180& 21003& 93067& 230848& & & & & \\
9& 1& 110& 4090& 55433& 333948& & & & & & \\
10& 1& 146& 7356& 137944& & & & & & & \\
\hline
\end{tabular}

\caption{ \texttt{( d)} unlabeled.
The number of directed graphs  allowing loops and multiedges
\cite[A138107]{EIS}.}
\end{table}

\begin{table}
\begin{tabular}{r|rrrrrrrrrrrrrrrrrr}
$E\backslash V$ & 1 & 2 & 3 & 4 & 5 & 6 & 7 &8 \\
\hline
0& 1& 1& 1& 1& 1& 1& 1& 1\\
1& 1& 4& 9& 16& 25& 36& 49& 64\\
2& 1& 10& 45& 136& 325& 666& 1225& 2080\\
3& 1& 20& 165& 816& 2925& 8436& 20825& 45760\\
4& 1& 35& 495& 3876& 20475& 82251& 270725& 766480\\
5& 1& 56& 1287& 15504& 118755& 658008& 2869685& 10424128\\
6& 1& 84& 3003& 54264& 593775& 4496388& 25827165& 119877472\\
7& 1& 120& 6435& 170544& 2629575& 26978328& 202927725& & & & \\
8& 1& 165& 12870& 490314& 10518300& 145008513& & & & & \\
9& 1& 220& 24310& 1307504& 38567100& & & & & & \\
10& 1& 286& 43758& 3268760& & & & & & & \\
\hline
\end{tabular}

\caption{
\texttt{( d)} vertex-labeled.
The number of labeled directed graphs  allowing loops and multiedges
\cite[A214398]{EIS}.}
\label{tab.fin}
\end{table}

\clearpage

\section{Connected Multigraphs up to 7 vertices}

\subsection{Algorithm}
The columns of the undirected connected multigraphs in Table \ref{tab.191646}
have rational ordinary generating functions. To compute them,
we first classify each multigraph by the number of edges and vertices
of the underlying simple graph---in as many ways as counted in Table \ref{tab.046742}---
and then distribute the edges of the multigraph over the edges of the
underlying graph using P\'olya's counting method to deal with the
symmetry of the simple graphs.

The process is illustrated in Section \ref{sec.4conm} for $V=4$ vertices.
Explicit intermediate results are tracked in the files \texttt{G.}$V$\texttt{.}$E$\texttt{.txt}
in the ancillary directory for $V=2$--$7$. Each of these files contains the contributing
underlying
simple graphs with $V$ vertices and $E$ edges. The file starts
with $V$ and $E$ printed in the first line. Then each graph is represented
by 
\begin{enumerate}
\item
a canonical adjacency matrix (binary, symmetric and traceless),
\item
the label as in Section \ref{sec.tag} followed by the multiplicity
of the graph as if one would create all vertex-labeled graphs by permuting
rows and columns (i.e. $V!$ divided by the order of the automorphism group),
\item
the cycle index multinomial.
This could also be derived from the table of symmetry groups in \cite[Vol. 1, Sec. 7, Table 8]{SteinbachVol4}.
\end{enumerate}
The minimum number of edges for connected simple graphs is $E\ge V-1$
(sparsest, trees on $V$ vertices), and the maximum number is $E\le \binom{V}{2}$
(complete graph on $V$ vertices). Summing over all multinomials
over the underlying graphs constitutes the generating
function by a finite sum of rational polynomials
(\cite[A001349]{EIS} terms).

\subsection{4 vertices}\label{sec.4conm}
The ordinary generating function 
for the number of connected multigraphs on 4 vertices
is derived
by adding the contributions of the 6 distinct geometries of the underlying
connected simple  graph.

We consider connected multigraphs with 4 vertices and
$E$ edges.
The multigraph thus has at least one (unoriented) edge
attached to edge vertex, so all degrees are $\ge 1$.
Loops are not allowed; the vertices are not labeled.

If all multiedges are replaced by a single edge,
the underlying simple graph has one of 6 shapes  \cite{GilbertCDM8}:

\begin{enumerate}
\item The linear chain with 3 edges.
\item A triangle with an edge to a lone vertex of degree 1 (4 edges).
\item The star graph with 3 edges.
\item The quadrangle (cycle of 4 edges).
\item A quadrangle with a single diagonal, total of 5 edges.
\item The complete graph on 4 vertices with 6 edges.
\end{enumerate}

\paragraph{Linear Chain.} 
We wish to distribute $E\ge 3$ edges over the 3 edges
of the simple graph of the linear chain. This can be done
by putting any number of $k\ge 1$ edges in the middle,
and distributing the remaining $n-k$ edges over the 
two edges connected to two endpoints.

Due to the left-right symmetry of the graph, the distribution
of the $n-k$ edges can only be done in $\lfloor (n-k)/2\rfloor$
ways. The total number of graphs of this kind with multiedges is
\begin{equation}
\sum_{k=1}^{E-2} \lfloor \frac{E-k}{2}\rfloor
=0, 0, 0, 1, 2, 4, 6, 9, 12, 16, 20, 25,\ldots (E\ge 0)
\end{equation}
with generating function \cite[A002620]{EIS}
\begin{equation}
g_1(x) = \sum_{E\ge 0} a_1(E) x^E = \frac{x^3}{(1+x)(1-x)^3}.
\end{equation}

The generating function is the product of $t_1(x)$ representing
the number of ways of placing $n$ vertices at the middle edge,
by the factor $x^2/[(1+x)(1-x)^2]$. The latter factor is
obtained by considering the symmetry of the cyclic group $C_2$
that swaps the edges that inhabit the first and last edges of the
underlying simple graph without generating a new graph. The cycle
index of the group is \cite{FreudensteinJM3}
\begin{equation}
Z(C_2) = (t_1^2+t_2^1)/2 ,
\end{equation}
where the associated generating functions are the number of
ways of placing $n$ edges without imposing symmetry on any
of them:
\begin{equation}
t_i(x) = \frac{x}{1-x} \mapsto 0,1,1,1,1,1\ldots, \quad i\ge 1.
\label{eq.tix}
\end{equation}
So the latter factor can be written as \cite[A004526]{EIS}
\begin{equation}
\frac{t_1(x)^2+t_2(x^2)}{2} = \frac{x^2}{(1+x)(1-x)^2}.
\label{eq.c2x}
\end{equation}

\paragraph{Triangle.} 

The contribution from the triangular graph is the number of ways
of placing $2\le k\le n-2$ edges on the edge to the lone vertex
and the triangle edge opposite to it, and then distributing the residual $n-k$
edges to the remaining two edges under the symmetry constraint of the group $C_2$
that swaps the other two edges:
\begin{equation}
\sum_{k=2}^{E-2} (k-1)\lfloor \frac{E-k}{2} \rfloor
= 0, 0, 0, 0, 1, 3, 7, 13, 22, 34, 50, 70, 95, 125, 161 (E\ge 0)
\end{equation}
See \cite[A002623]{EIS}
\begin{equation}
g_2(x) = \frac{x^4}{(1+x)(1-x)^4}.
\label{eq.g2x}
\end{equation}
This generating function is the product of $x^2/(1-x)^2$---contribution
of two cycles of length 1, fixed points under the symmetry---by
$x^2/[(1+x)(1-x)^2]$, where again the latter is \eqref{eq.c2x},
the number of ways of distributing $E$ edges symmetrically 
over two edges of the simple graph.

\paragraph{Star Graph.} 

The contribution from the star graph is the number of ways
of partitioning $E$ into 3 positive integers \cite[A069905]{EIS},
\begin{equation}
\mapsto 0, 0, 0, 1, 1, 2, 3, 4, 5, 7, 8, 10, 12, 14, 16 (E\ge 0),
\end{equation}
\begin{equation}
g_3(x) = \frac{x^3}{(1+x)(1-x)^2(1-x^3)}.
\label{eq.g3x}
\end{equation}
Alternatively this expression is obtained if we consider
the symmetry group of order 6 of the underlying simple graph,
which can be generated by (i) the group $C_3$ of the triangle
combined with (ii) the mirror symmetry along a diagonal.
\begin{verbatim}
with(group):
g := permgroup(3, {[[1, 2, 3]], [[2, 3]]}) ;
for i in elements(g) do
    print(i) ;
end do;
\end{verbatim}
The cycle index obtained with this Maple code is
\cite[p 57]{Polya1983}
\begin{equation}
Z(S_3) = \frac{t_1^3 +3t_1t_2+2t_3^1}{6}.
\end{equation}
Insertion of \eqref{eq.tix} gives \eqref{eq.g3x}.

\paragraph{Square.} 
The symmetry group of the square is the Dihedral Group of order 8 which
essentially is generated by rotation by 90 degrees or flips
along the horizontal or vertical axes or diagonals.
\begin{verbatim}
with(group):
g := permgroup(4, {[[1, 2, 3, 4]], [[1, 3]]}) ;
for i in elements(g) do
    print(i) ;
end do;
\end{verbatim}
The cycle index is \cite[Fig 3]{TuckerMM47}\cite{FreudensteinJM3}
\begin{equation}
Z(D_8) = \frac{t_1^4+2t_4^1+2t_1^2t_2^1+3t_2^2}{8}
\end{equation}
The enumeration theorem turns this into the generating function
\cite[A005232]{EIS}
\begin{equation}
g_4(x) = \frac{x^4(1-x+x^2)}{(1+x^2)(1+x)^2(1-x)^4}.
\end{equation}
\begin{equation}
\mapsto 0, 0, 0, 0, 1, 1, 3, 4, 8, 10, 16, 20, 29, 35, 47 (E\ge 0).
\end{equation}

\paragraph{Square with Diagonal.} 
In the square with a diagonal edge, the diagonal stays inert under
the symmetry operations, and contributes a factor $t_1(x)$ to the
generating function. The symmetry group of the four other edges
allows a flip along any of the two diagonals and generates a symmetry group
of order 4:
\begin{verbatim}
with(group):
g := permgroup(4,{[[2, 4], [1, 3]], [[1, 4], [2, 3]]}) ;
for i in elements(g) do
    print(i) ;
end do;
\end{verbatim}
The cycle index is
\begin{equation}
Z(C_2\times C_2) = \frac{t_1^4+3t_2^2}{4}.
\end{equation}
Insertion of \eqref{eq.tix} into the enumeration theorem
yields $x^4(1-x+x^2)/[(1+x)^2(1-x)^4]$ \cite[A053307]{EIS}, and convolved with the inert
factor
\begin{equation}
g_5(x) = \frac{x^5(1-x+x^2)}{(1+x)^2(1-x)^5}.
\end{equation}
This expands to
\begin{equation}
\mapsto 0, 0, 0, 0, 0, 1, 2, 6, 11, 22, 36, 60, 90, 135, 190, 266, 357, 476 (E\ge 0)
\end{equation}

\paragraph{Complete Graph.} 
The cycle index of the complete graph $K_4$ 
is
\cite{FreudensteinJM3},
\begin{equation}
Z(S_4)=\frac{t_1^6+9t_1^2t_2^2+8t_3^2+6t_2t_4}{24}.
\end{equation}
Insertion of \eqref{eq.tix} into the enumeration theorem yields
\begin{equation}
g_6(x) = \frac{x^6(1-x+x^2+x^4+x^6-x^7+x^8)}{(1-x)^6(1+x)^2(1+x^2)(1+x+x^2)^2},
\end{equation}
with \cite[A003082]{EIS}
\begin{equation}
\mapsto 0, 0, 0, 0, 0, 0, 1, 1, 3, 6, 11, 18, 32, 48, 75, 111, 160, 224, 313 (E\ge 0)
\end{equation}

\paragraph{Sum.}
The generating function contributed by the 6 underlying simple graphs is
\begin{equation}
\sum_{i=1}^6 g_i(x)
=
\frac{x^3(-x^{10}+x^9+2x^7-x^6+x^5-3x^4+x^2+x+2)}{(x-1)^6(1+x)^2(1+x^2)(1+x+x^2)^2}
,
\end{equation}
which expands to \cite[A290778]{EIS}
\begin{equation}
\mapsto
0, 0, 0, 2, 5, 11, 22, 37, 61, 95, 141, 203, 288, 393, 531, 704, 918, 1180, 1504
(E\ge 0, V=4)
\end{equation}

\subsection{Up to 7 vertices}

The generating function for the number of connected loopless multigraphs
on 2 vertices  is
\begin{equation}
\frac{x}{1-x}
\mapsto 0,1,1,1,1,1, (E\ge 0, V=2)
.
\end{equation}

The generating function for the number of connected loopless multigraphs
on 3 vertices is \cite[A253186]{EIS}
\begin{multline}
{\frac { \left( {x}^{3}-x-1 \right) {x}^{2}}{ \left( -1+x \right) ^{3}
 \left( x+1 \right)  \left( {x}^{2}+x+1 \right) }}
\\
\mapsto
0, 0, 1, 2, 3, 4, 6, 7, 9, 11, 13, 15, 18, 20, 23 (E\ge 0, V=3)
.
\end{multline}

On 5 vertices
\begin{multline}
{\frac {{x}^{4} 
p_5(x)
}{ \left( -1+x \right) ^{10} \left( {x}^{2}+x+1 \right) ^{3}
 \left( x+1 \right) ^{4} \left( {x}^{2}-x+1 \right)  \left( {x}^{4}+{x}
^{3}+{x}^{2}+x+1 \right) ^{2} \left( {x}^{2}+1 \right) ^{2}}}
\\
\mapsto
0, 0, 0, 0, 3, 11, 34, 85, 193, 396, 771, 1411, 2490, 4221 (E\ge 0, V=5),
\end{multline}
where
\begin{multline}
p_5 \equiv 
3+5\,x+12\,{x}^{2}+17\,{x}^{3}+26\,{x}^{4}+27\,{
x}^{5}+35\,{x}^{6}+28\,{x}^{7}+38\,{x}^{8}+30\,{x}^{9}+39\,{x}^{10}\\
+37 \,{x}^{11}
+34\,{x}^{12}+24\,{x}^{13}+15\,{x}^{14}+3\,{x}^{15}
-7\,{x}^{
16}-9\,{x}^{17}+4\,{x}^{20}\\
+5\,{x}^{22}+3\,{x}^{23}-8\,{x}^{18}-{x}^{19
}+6\,{x}^{21}-2\,{x}^{24}-2\,{x}^{25}-2\,{x}^{26}-{x}^{27}+{x}^{29}
.
\end{multline}
On 6 vertices
\begin{multline}
-{\frac {{x}^{5} p_6
}{ \left(
-1+x \right) ^{15} \left( x+1 \right) ^{6} \left( {x}^{2}+1 \right) ^{3
} \left( {x}^{2}+x+1 \right) ^{5} \left( {x}^{2}-x+1 \right) ^{2}
 \left( {x}^{4}+{x}^{3}+{x}^{2}+x+1 \right) ^{3}}}
\\
\mapsto 0, 0, 0, 0, 0, 6, 29, 110, 348, 969, 2445, 5746, 12736, 26843, 54256, 105669 (E\ge 0, V=6),
\end{multline}
where
\begin{multline}
p_6 \equiv
 6
+11\,x
+35\,{x}^{2}
+70\,{x}^{3}
+134\,{x}^{4}
+ 217\,{x}^{5}
+348\,{x}^{6}
+533\,{x}^{7}
+726\,{x}^{8}
+1038\,{x}^{9}
+1290 \,{x}^{10} 
\\
+1629\,{x}^{11}
+1810\,{x}^{12}
+2040\,{x}^{13}
+1976\,{x}^{14}
+ 1984\,{x}^{15}
+1696\,{x}^{16}
+1542\,{x}^{17}
+1206\,{x}^{18}
\\
+1050\,{x}^{19} 
+787\,{x}^{20}
+636\,{x}^{21}
+474\,{x}^{22 }
+273\,{x}^{23}
+169\,{x}^{ 24}
-11\,{x}^{25}
-31\,{x}^{26}
-97\,{x}^{27}
-44\,{x}^{28}
\\
-8\,{x}^{29}
+33 \,{x}^{30} 
+63\,{x}^{31}
+32\,{x}^{32}
+38\,{x}^{33}
-17\,{x}^{34}
-14\,{x}^ {35}
-31\,{x}^{36}
-8\,{x}^{37}
-5\,{x}^{38}
+8\,{x}^{39}
\\
+11\,{x}^{40}
+4\,{ x}^{41}
+3\,{x}^{42}-4\,{x}^{43}-3\,{x}^{45}
+{x}^{47}
.
\end{multline}
On 7 vertices
\begin{multline}
\frac {{x}^{6}p_7(x)
}{ \left( -1+x \right) ^{21} \left( {x
}^{4}+{x}^{3}+{x}^{2}+x+1 \right) ^{4} \left( {x}^{2}+x+1 \right) ^{7}
 \left( x+1 \right) ^{9} \left( {x}^{2}+1 \right) ^{4} \left( {x}^{2}-x
+1 \right) ^{3}
}
\\
\times
\frac{1}{
\left( {x}^{4}-{x}^{2}+1 \right)  \left( {x}^{4}-{x}^{3
}+{x}^{2}-x+1 \right)  \left( {x}^{6}+{x}^{5}+{x}^{4}+{x}^{3}+{x}^{2}+x
+1 \right) ^{3}}
\\
\mapsto
0, 0, 0, 0, 0, 0, 11, 70, 339, 1318, 4457, 13572, 38201, 100622, 251078,
    \ldots (E\ge 0, V=7),
\end{multline}
where
\begin{multline}
p_7(x)\equiv
-11
-48\,x
-188\,{x}^{2}
-570\,{x}^{3}
-1526\,{x}^{4 }
-3675\,{x}^{5}
-8284\,{x}^{6}
-17431\,{x}^{7}
-35005\,{x}^{8}
\\
-66742\,{x}^ {9}
-121908\,{x}^{10}
-213342\,{x}^{11}
-359515\,{x}^{12}
-583522\,{x}^{13} 
-916091\,{x}^{14}
-1391716\,{x}^{15}
\\
-2051981\,{x}^{16}
-2938963\,{x}^{17} 
-4097420\,{x}^{18}
-5564508\,{x}^{19}
-7373793\,{x}^{20}
-9539279\,{x}^{21 }
\\
-12063528\,{x}^{22}
- 14919997\,{x}^{23}
-18064473\,{x}^{24}
-21418776\,{x}^{25}
-24890827\,{x}^{26}
-28355984\,{x }^{27}
\\
-31688266\,{x}^{28}
-34742272\,{x}^{29}
-37387611\,{x}^{30}
- 39493274\,{x}^{31}
-40963946\,{x}^{32}
-41717383\,{x}^{33}
\\
-41723196\,{x}^ {34}
-40973187\,{x}^{35}
-39511812\,{x}^{36}
-37405689\,{x}^{37}
-34764514 \,{x}^{38}
-31705308\,{x}^{39}
\\
-28372262\,{x}^{40}
-24898844\,{x}^{41}
- 21423490\,{x}^{42}
-18060699\,{x}^{43}
-14913079\,{ x}^{44}
-12050303\,{x}^{45}
\\
-9525196\,{x}^{46} 
-7357519\,{x}^{ 47}
-5550815\,{x}^{48}
-4085547\,{x}^{49}
-2932089\,{x}^{50}
-2048825\,{x}^{51}
\\
-1393454\,{x} ^{52}
-920594\,{x}^{53}
-590477\,{x}^{54}
-366935\,{x }^{55}
-220705\,{x}^{56}
-128024\,{x}^{57}
-71511\,{x}^{58 }
\\
-37993\,{x}^{59}
-18932\,{x}^{60}
-8318\,{x}^{61}
-2668\,{x}^{62}
+247\,{x}^{63}
+1501\,{x}^{64}
+1827\,{x}^{65}
+1523\,{x}^ {66}
\\
+980\,{x}^{67}
+357\,{x}^{68}
-99\,{x}^{69}
-369\,{x}^{70}
-387\,{x}^{71}
-247\,{x}^{ 72}
-23\,{x}^{73}
+152\,{x}^{74}
\\
+230\,{x}^{75}
+205\,{x}^{76}
+118\,{x}^{77 }
+15\,{x}^{78}
-61\,{x}^{79}
-88\,{x}^{80}
-74\,{x}^{81}
-33\,{x}^{82}
+3\,{x}^{83}
+26\,{x}^{84}
\\
+28\,{ x}^{85}
+19\,{x}^{86}
+5\,{x}^{87}
-4\,{x}^{88}
-7\,{x}^{89}
-5\,{x}^{90}
-{x}^{91}
+{x}^ {92}
+{x}^{93}
.
\end{multline}

As a cross-check on these numbers we note that using a weight $t_i(x)=x$ instead
of \eqref{eq.tix} just counts the underlying simple graphs; it computes the generating functions
down columns of Table \ref{tab.046742}.

\section{Multisets of Connected Graphs}

If a table of connected graphs as a function of vertex count and edge count
is known, the Multiset Transformation generates tables of disconnected graphs
with fixed number of components. The calculation involves creating an
intermediate multiset of edge-vertex pairs of the components, and looking
up a product of multiset coefficients as a function of the number of connected
graphs that support the pairs.   The technique is demonstrated for
simple (undirected, unlabeled, loopless) graphs and for  undirected, unlabeled loopless
graphs allowing multiedges.

\subsection{The Multiset Coefficient}

The concept of the multiset is based on the concept of the set (a collection
of objectes, only one object of a given type), but allows to put more than
one object of a type into the collection \cite{SinghNSJM37}.

\begin{defn}
A Multiset is a collection of objects with some individual
count (object of type $i$ appearing $f_i$ times in the collection).
The objects have no order in the collection.
\end{defn}

The number of ways of assembling a multiset with $m$ objects plugged from a
set of $n$ different objects is a variant of Pascal's triangle
of binomial coefficients,
\begin{equation}
P(n,m)=\binom{n+m-1}{m}
.
\label{eq.P}
\end{equation}
The equation may be illustrated for small orders $m$:
\begin{itemize}
\item
If there is only $n=1$ type of objects, the multiset has only one choice: it contains $m$
replicates of the unique object. $P(1,m)=1$.
\item
If the multiset contains $m=1$ object, it contains
one object of $n$ candidates. 
$P(n,m)=n$.
\item
If the multiset contains $m=2$ objects, it
contains either the same type of object twice ($n$ choices), or two different
objects ($\binom{n}{2}$ choices), so $P(n,2)=n+\binom{n}{2} = \binom{n+1}{2}$.
\item
If the multiset contains $m=3$ objects, 
it either contains the same type of object thrice ($n$ choices), 
or one type of object once and another type of object twice ($n(n-1)$ choices),
or three different types of objects ($\binom{n}{3}$ choices);
so $P(n,3)=n
+n(n-1)
+\binom{n}{3}
=\binom{n+2}{3}$.
\item
If the multiset contains $m=4$ objects, we consider all five
partitions of $m$, namely ${4^1},{1^13^1},{2^2},{1^22^1},{1^4}$:
it either contains the same type of object 4 times ($n$ choices), 
or one type of object once and another type of object thrice ($n(n-1)$ choices),
or two pairs of objects ( $\binom{n}{2}$ choices),
or two different objects and one pair of objects ( $\binom{n}{2}(n-2)$ choices),
or four different types of objects ($\binom{n}{4}$ choices);
so $P(n,4)=n +n(n-1) +\binom{n}{2}+\binom{n}{2}(n-2)+\binom{n}{4} =\binom{n+3}{4}$.
\end{itemize}

\begin{proof}
Formula \eqref{eq.P} can be rephrased with \cite[(3.8)]{RoyAMM94}
(setting $r=m-1$, $k=j$, $n=1$ there) or with \cite[(1.11)]{Sane}:
\begin{equation}
\binom{n+m-1}{m}
=
\sum_{j=1}^m \binom{n}{j} \binom{m-1}{j-1}
.
\end{equation}
The first factor on the right hand side indicates 
that in a first step one can create a set of $j$ distinct objects out of $n$ in $\binom{n}{j}$ ways.
Consider that set sorted by some lexicographic order. Then the factor $\binom{m-1}{j-1}$ counts
in how many ways one can insert separators in the multiset of the same lexicographic ordering to
select switch-over from one type of object to the next one.
\end{proof}
\begin{exa}
To create multisets of $m=4$ objects given $n$ distinct objects, selecting $j=1$ type of object
gives the multiset $o_1o_1o_1o_1$ (no separator), $\binom{3}{0}=1$ choices;
selecting $j=2$ types of objects gives
$o_1o_1o_1|o_2$ or $o_1o_1|o_2o_2$ or $o_1|o_2o_2o_2$ with $\binom{3}{1}=3$ positions of the separator;
selecting $j=3$ types of objects gives
$o_1o_1|o_2|o_3$ or $o_1|o_2o_2|o_3$ or $o_1|o_2|o_3o_3$ with $\binom{3}{2}=3$ positions of the $2$ separators;
or selecting $j=4$ types of objects gives
$o_1|o_2|o_3|o_4$ with $\binom{3}{3}=1$ positions of the $3$ separators.
\end{exa}

\subsection{Multiset Transform of An Integer Sequence}
The Multiset Transform deals with the question: if the objects of
type $i$ have some weight (expressed as a positive integer), in how
many ways can we assemble a multiset of the objects with some prescribed
total weight? The total weight is the usual arithmetic sum of the weights of the objects.
\begin{exa}
(Money Exchange Problem)
In how different ways can you combine coins (weight 5 for type nickels, weight 10
for type dimes and weight 25 for type quarter\ldots) for a wallet worth 200
(2 dollars)?
\end{exa}
\begin{exa}
In how different ways can you fill a bag of 10 kg (weight 100) with
apples of 100 g (weight 1) and oranges of 200 g (weight 2)?
\end{exa}

The Multiset Transform computes the number $T_{n,k}$ of multisets containing
$k$ objects, drawn from a set of objects of which there are $T_{n,1}$ of some
additive weight $n$ \cite{Flajolet}.  Given $T_{0,1}=1$ and an integer sequence $T_{n,1}$
for the number of objects in the weight class $n$, the $T$ of total
weight $n$ are calculated
recursively by
\begin{equation}
T_{n,k}=\sum_{f_1n_1+f_2n_2+\cdots +f_kn_k=n} \prod_{i=1}^k P(T_{n_i,1},f_i),
\label{eq.mrec}
\end{equation}
where the sum is over the partitions of $n$ into parts $n_i$ which 
occur with frequencies $f_i$. 

\begin{remark}
The row sums $\sum_{k=1}^n T_{n,k}$ are obtained by the Euler
Transform of the sequence $T_{n,1}$ \cite[(25)]{Flajolet}.
\end{remark}

\begin{exa}
Given two types of nickels (2 types of weight 5, tin and copper), one type
of dime (weight 10), and one type of quarter (weight 25), the integer
sequence $T_{n,1}$ is $(1),0,0,0,0,2,0,0,0,0,1,0,0,0,0,0,0,0,0,0,0,0,0,0,0,25$.
\end{exa}

\begin{exa}
Given three types of apples (brown, yellow and red, each of weight 1), one type
of banana (weight 2), and one type melon (weight 4), the integer
sequence is $(1),3,1,0,1,0,0,0\ldots$. The Multiset Transform generates
the triangular table

\begin{tabular}{r|rrrrrrrrr}
$n\backslash k$ & 1 & 2 & 3 & 4 & 5 & 6 \\
\hline
1& 3\\
2& 1&6\\
3& 0&3&10\\
4& 1&1&6&15\\
5& 0&3&3&10&21\\
6& 0&1&7&6&15&28\\
7& 0&0&3&13&10&21&36\\
8& 0&1&1&7&21&15&28&45\\
\end{tabular}

The row sums in the table are 3, 7, 13, 23, 37, 57, 83, 118\ldots
There are $T_{3,2}=3$ ways of generating a weight of 3 with two 
objects (a banana and any of the three types of apples).
There are $T_{3,3}=10$ ways of generating a weight of 3 with three
objects [three apples (bbb), (yyy), (rrr), (byy), (brr), (bby), (yrr), (bbr), (yyr), (bry)].
There is $T_{6,2}=1$ way to generate a weight of 6 with two objects (a banana and a melon).
\end{exa}

\begin{exa}
If there is one type of object of each weight, $T_{n,1}=1$, the Multiset Transform generates
the partition numbers \cite[A008284]{EIS}, and the row sums are the partition numbers \cite[A000041]{EIS}.
\end{exa}

\subsection{Graphs specified by number of components}\label{sec.grexa}

If the sequence $T_{n,1}$ enumerates connected graphs of type $n$ (where
the weight $n$ is either the vertex count or the edge count), one fundamental
way of generating a multiset is putting $k$ of them side by side
and considering them a graph with $k$
components. Vertex or edge number are additive, as required.

\begin{exa}
If $T_{n,1}$ denotes connected graphs with $n$ nodes \cite[A001349]{EIS},
the Multiset Transform counts graphs with $k$ components \cite[A201922]{EIS}.
\end{exa}

\begin{exa}
If $T_{n,1}$ denotes connected graphs with $n$ edges \cite[A002905]{EIS},
the Multiset Transform counts graphs with $n$ edges and $k$ components \cite[A275421]{EIS}.
\end{exa}

\begin{exa}
If $T_{n,1}$ denotes trees with $n$ nodes \cite[A000055]{EIS},
the Multiset Transform counts forests with $k$ trees \cite[A095133]{EIS}.
\end{exa}

\begin{exa}
If $T_{n,1}$ denotes rooted trees with $n$ nodes \cite[A000081]{EIS},
the Multiset Transform counts rooted forests with $k$ trees \cite[A033185]{EIS}.
\end{exa}

\begin{exa}
If $T_{n,1}$ denotes connected regular graphs with $n$ nodes \cite[A005177]{EIS},
the Multiset Transform counts regular graphs with $k$ components \cite[A275420]{EIS}.
In the case of cubic graphs the transform pair is \cite[A002851]{EIS} and \cite[A275744]{EIS}.
\end{exa}

\section{Graphs specified by number of edges, vertices and components}

\subsection{Union of connected graphs}
Let $G(E,V,k)$ be the number of graphs with $E$ edges, $V$ vertices
and $k$ components. $G(E,V,1)$ is the number of connected graphs
with $E$ edges and $V$ vertices. The other properties like whether the
graphs are labeled, may contain loops or multiedges, are
not classified here, but assumed to be fixed while composing
graphs with $k$ components from connected graphs. (A multiset of labeled connected
graphs is a disconnected labeled graph;
a multiset of connected oriented graphs is a disconnected oriented graph; and
so on.)
The unified
graph is the multiset union of connected graphs $\mathcal G_1$
which individually have $e_i$ edges and $v_i$ vertices:
\begin{equation}
\mathcal G_k(E,V) = \bigcup_{i=1}^k \mathcal G_1(e_i,v_i),
\end{equation}
where both the number of edges and the number of vertices are
additive:
\begin{equation}
E=\sum_i^k e_i;\quad V=\sum_i^k v_i.
\label{eq.add}
\end{equation}
Summation over a set of the variables creates marginal sums:
$G(.,V,k)=\sum_{E\ge 0}G(E,V,k)$ are the graphs with $V$ vertices
and $k$ components. $G(E,,k)=\sum_{V\ge 1} G(E,V,k)$ is the number
of graphs with $E$ edges and $k$ components. $G(E,V,.)=\sum_{k\ge 1} G(E,V,k)$
is the number of graphs with $E$ edges and $V$ vertices.

\subsection{Correlated Multiset Transforms}
The particular case we explore here is that constructing
a disconnected graph from connected components means building a multiset
of connected graphs, where the number of edges and \emph{also} the number of
vertices are such an additive weight. 

The examples of Section \ref{sec.grexa}
illustrated how
$G(E,,k)$ is the Multiset Transform of $G(E,,1)$ 
\cite[A076864,A275421,A191970]{EIS} and
$G(,V,k)$ is the Multiset Transform of $G(,V,1)$ \cite[A054924,A275420,A281446]{EIS}. 
The aim of this paper is to demonstrate
a similar technique for $G(E,V,k)$ assuming $G(E,V,1)$ is known.

Each graph which is a component contributing to $G(E,V,k)$ has
a specific pair $(e_i,v_i)$ of edge and vertex count; the union of these
graphs is a multiset of such pairs---which means in the multiset of graphs contributing
to $G(E,V,k)$, each pair may occur more than once, and each pair may represent
(in the sense of the weights above) more than one graph because there
may be more than one distinct connected graph for one pair of $E$ and $V$.

The calculation starts by constructing all weak compositions of $E$
into $k$ parts $e_i$, and all weak compositions of $V$ into $k$
parts $v_i$ of pairs $(e_i,v_i)$ compatible with the requirement \eqref{eq.add}.
This defines a two-dimensional $\binom{E+k-1}{k-1}\times \binom{V+k-1}{k}$
outer product matrix  with multisets \cite{WiederPAM2}.
\begin{exa}\label{exa.EV}
If $E=2$ and $V=3$ and $k=3$, the compositions are
$2=$$2+0+0=$$0+2+0=$$0+0+2=$$1+1+0=$$1+0+1=$$0+1+1$,
$3=$$3+0+0=$$0+3+0=$$0+0+3=$$2+1+0=$$1+2+0=$$2+0+1=$$2+1+0=\ldots$, 
and the matrix contains multisets with $k$ pairs:
\begin{tabular}{r|rrrrrrrrrrrrrrrrr}
$\sum e_i\backslash \sum v_i$ &       3+0+0 &          0+3+0 & & 2+1+0 & $\ldots$ & 1+1+1 \\
\hline
2+0+0 & (2,3)(0,0)(0,0) & (2,0)(0,3)(0,0) & \ldots & & & (2,1)(0,1)(0,1) \\
0+2+0 & (0,3)(2,0)(0,0) & (0,0)(2,3)(0,0) & \ldots & & & (0,1)(2,1)(0,1) \\
0+0+2 & \ldots \\
1+1+0 & \ldots \\
0+1+1 & (0,3)(1,0)(1,0) & (0,0)(1,3)(1,0) & \ldots & & & (0,1)(1,1)(1,1) \\
& $\ldots$ \\
\hline
\end{tabular}
\end{exa}
Each element of the matrix is a multiset $(e_1,v_1)(e_2,v_2)\cdots (e_k,v_k)$
of pairs obtained by interleaving the $e_i$ and $v_i$
components of the compositions. At that point we realize that 
\begin{enumerate}
\item
if any of the $v_i$ is zero, the method of selecting such a 
null-graph into the disconnected graph would not fulfill the
requirement of being $k$-connected. So actually only the
compositions (not the weak compositions) of $V$ need to be considered
as table columns.
\item
because the decompositions $(e_1,v_1),(e_2,v_2)\cdots$
are candidates for multiset compositions, the order of
the pairs does not matter. In the example, $(2,1)(0,1)(0,1)$
and $(0,1),(2,1)(0,1)$ are the same scheme of selecting connected
graphs. So we may sort for example each $k$-set by the $e_i$
member of the pair without loss of samples, which means,
we may build the table by just considering weak partitions
(not all weak compositions) of $E$ into parts $e_i$ as table rows.
\end{enumerate}
\begin{exa}
Continuing Example \ref{exa.EV} above, the
table reduces to 

\begin{tabular}{r|r}
$\sum e_i\backslash \sum v_i$ &       1+1+1 \\
\hline
2+0+0 & (2,1)(0,1)(0,1) \\
1+1+0 & (1,1)(1,1)(0,1) \\
\hline
\end{tabular}

The $(2,1)(0,1)(0,1)$ entry indicates to take one connected
graph with 2 edges and one vertex (obviously a double-loop) and
two graphs without edges and one vertex (single points). The $(1,1)(1,1)(0,1)$
entry indicates to take 2 graphs with an edge and a vertex (2 points,
each with a loop) and a graph without edges and a single vertex (a single point).
\end{exa}

The number of $k$-compositions of graphs then is the
sum over all unique remaining multisets in the table, akin
to Equation \eqref{eq.mrec}:
\begin{equation}
G(E,V,k) = \sum_{(e_i,v_i)^{f_i}}
\prod_i P(G(e_i,v_i,1),f_i)
\end{equation}
where $f_i$ is the frequency (number of occurrences) of the pair
$(e_i,v_i)$ in the multiset.

\subsection{Simple Graphs}
The most important example considers undirected, unlabeled graphs without 
multiedges or loops. Table \ref{tab.046742} shows the base
information $G(E,V,1)$ which is fed into the formula
to compute the tables $G(E,V,k)$  \ref{tab.sk2}--\ref{tab.sk5}
for $k\ge 1$.
(These tables differ from Steinbach's tables \cite[Vol. 4, Table 1.2a]{SteinbachVol4}
because our components may be/have isolated vertices.)

\begin{table}
\begin{tabular}{r|rrrrrrrrrrrrrrrrrrr}
$E\backslash V$ & 1 & 2 & 3 & 4 & 5 & 6 & 7 &8 &9 &10 &11 & 12 & 13 &14 \\
\hline
0 &0 &1\\
1 &0 &0 &1\\
2 &0 &0 &0 &2\\
3 &0 &0 &0 &1 &3\\
4 &0 &0 &0 &0 &3 &6\\
5 &0 &0 &0 &0 &1 &8 &11\\
6 &0 &0 &0 &0 &1 &7 &22 &23\\
7 &0 &0 &0 &0 &0 &5 &27 &58 &46\\
8 &0 &0 &0 &0 &0 &2 &28 &101 &157 &99\\
9 &0 &0 &0 &0 &0 &1 &23 &142 &358 &426 &216\\
10 &0 &0 &0 &0 &0 &1 &15 &161 &660 &1233 &1166 &488\\
11 &0 &0 &0 &0 &0 &0 &10 &156 &1010 &2873 &4163 &3206 &1121\\
12 &0 &0 &0 &0 &0 &0 &5 &138 &1356 &5705 &11987 &13847 &8892 &2644\\
13 &0 &0 &0 &0 &0 &0 &2 &101 &1613 &9985 &29652 &48071 &45505 &24743 \\

\end{tabular}
\caption{
$G(E,V,2)$. Simple graphs with 2 components and a total of $E$ edges and $V$ vertices.
See \cite[A274934]{EIS} for column sums, \cite[A274937]{EIS} for the diagonal.
}
\label{tab.sk2}
\end{table}

\begin{table}
\begin{tabular}{r|rrrrrrrrrrrrrrrrrrr}
$E\backslash V$ & 1 & 2 & 3 & 4 & 5 & 6 & 7 &8 &9 &10 &11 & 12 & 13 &14 \\
\hline
0 &0 &0 &1\\
1 &0 &0 &0 &1\\
2 &0 &0 &0 &0 &2\\
3 &0 &0 &0 &0 &1 &4\\
4 &0 &0 &0 &0 &0 &3 &7\\
5 &0 &0 &0 &0 &0 &1 &9 &14\\
6 &0 &0 &0 &0 &0 &1 &7 &25 &29\\
7 &0 &0 &0 &0 &0 &0 &5 &29 &68 &60\\
8 &0 &0 &0 &0 &0 &0 &2 &29 &110 &186 &128\\
9 &0 &0 &0 &0 &0 &0 &1 &23 &149 &397 &509 &284\\
10 &0 &0 &0 &0 &0 &0 &1 &15 &164 &699 &1377 &1399 &636\\
11 &0 &0 &0 &0 &0 &0 &0 &10 &157 &1041 &3070 &4685 &3857 &1467\\
12 &0 &0 &0 &0 &0 &0 &0 &5 &139 &1375 &5919 &12899 &15646 &10706 \\
13 &0 &0 &0 &0 &0 &0 &0 &2 &101 &1625 &10183 &30980 &52024 &51622 \\

\end{tabular}
\caption{
$G(E,V,3)$. Simple graphs with 3 components and a total of $E$ edges and $V$ vertices.
Column sums are in column 3 of \cite[A201922]{EIS} $\mapsto 1,1,3,9,32,154,1065,12513,276114,12021725\ldots$.
}
\end{table}

\begin{table}
\begin{tabular}{r|rrrrrrrrrrrrrrrrrrr}
$E\backslash V$ & 1 & 2 & 3 & 4 & 5 & 6 & 7 &8 &9 &10 &11 & 12 & 13 &14 & 15\\
\hline
0 &0 &0 &0 &1\\
1 &0 &0 &0 &0 &1\\
2 &0 &0 &0 &0 &0 &2\\
3 &0 &0 &0 &0 &0 &1 &4\\
4 &0 &0 &0 &0 &0 &0 &3 &8\\
5 &0 &0 &0 &0 &0 &0 &1 &9 &15\\
6 &0 &0 &0 &0 &0 &0 &1 &7 &26 &32\\
7 &0 &0 &0 &0 &0 &0 &0 &5 &29 &71 &66\\
8 &0 &0 &0 &0 &0 &0 &0 &2 &29 &112 &196 &143\\
9 &0 &0 &0 &0 &0 &0 &0 &1 &23 &150 &406 &539 &315\\
10 &0 &0 &0 &0 &0 &0 &0 &1 &15 &164 &706 &1417 &1486 &710\\
11 &0 &0 &0 &0 &0 &0 &0 &0 &10 &157 &1044 &3110 &4834 &4105 &1631\\
12 &0 &0 &0 &0 &0 &0 &0 &0 &5 &139 &1376 &5951 &13102 &16193 &11408 \\
13 &0 &0 &0 &0 &0 &0 &0 &0 &2 &101 &1626 &10202 &31198 &52966 &53519 \\

\end{tabular}
\caption{
$G(E,V,4)$. Simple graphs with 4 components and a total of $E$ edges and $V$ vertices.
}
\label{tab.sk4}
\end{table}

\begin{table}
\begin{tabular}{r|rrrrrrrrrrrrrrrrrrr}
$E\backslash V$ & 1 & 2 & 3 & 4 & 5 & 6 & 7 &8 &9 &10 &11 & 12 & 13 &14 & 15\\
\hline
0 &0 &0 &0 &0 &1\\
1 &0 &0 &0 &0 &0 &1\\
2 &0 &0 &0 &0 &0 &0 &2\\
3 &0 &0 &0 &0 &0 &0 &1 &4\\
4 &0 &0 &0 &0 &0 &0 &0 &3 &8\\
5 &0 &0 &0 &0 &0 &0 &0 &1 &9 &16\\
6 &0 &0 &0 &0 &0 &0 &0 &1 &7 &26 &33\\
7 &0 &0 &0 &0 &0 &0 &0 &0 &5 &29 &72 &69\\
8 &0 &0 &0 &0 &0 &0 &0 &0 &2 &29 &112 &199 &149\\
9 &0 &0 &0 &0 &0 &0 &0 &0 &1 &23 &150 &408 &549 &330\\
10 &0 &0 &0 &0 &0 &0 &0 &0 &1 &15 &164 &707 &1426 &1516 &742\\
11 &0 &0 &0 &0 &0 &0 &0 &0 &0 &10 &157 &1044 &3117 &4874 &4193 \\
12 &0 &0 &0 &0 &0 &0 &0 &0 &0 &5 &139 &1376 &5954 &13142 &16343 \\
13 &0 &0 &0 &0 &0 &0 &0 &0 &0 &2 &101 &1626 &10203 &31230 &53170 \\

\end{tabular}
\caption{
$G(E,V,5)$. Simple graphs with 5 components and a total of $E$ edges and $V$ vertices.
}
\label{tab.sk5}
\end{table}

The arithmetic sum over these tables $k\ge 1$ yields Table \ref{tab.008406}.

\clearpage
\subsection{Loopless connected Multigraphs}
Another application of the algorithm is to construct
\cite[A192517]{EIS} from \cite[A191646]{EIS}.
The base information of $G(E,V,k=1)$ is this time
in Table \ref{tab.191646}, and the Multiset Transformations
creates Tables \ref{tab.mk2}--\ref{tab.mk5} and so on ($k\ge 2$).
The sum over all $k\ge 1$ converges to Table \ref{tab.192517}.

\begin{table}
\begin{tabular}{r|rrrrrrrrrrrrrrrrrrrr}
$E\backslash V$ & 1 & 2 & 3 & 4 & 5 & 6 & 7 &8 &9 &10 &11 & 12 \\
\hline
0 &0 &1\\
1 &0 &0 &1\\
2 &0 &0 &1 &2\\
3 &0 &0 &1 &3 &3\\
4 &0 &0 &1 &5 &8 &6\\
5 &0 &0 &1 &6 &17 &20 &11\\
6 &0 &0 &1 &9 &32 &58 &52 &23\\
7 &0 &0 &1 &10 &53 &135 &185 &132 &46\\
8 &0 &0 &1 &13 &84 &290 &548 &586 &344 &99\\
9 &0 &0 &1 &15 &127 &565 &1441 &2108 &1829 &900 &216\\
10 &0 &0 &1 &18 &184 &1055 &3456 &6696 &7884 &5680 &2834 &488\\
11 &0 &0 &1 &20 &259 &1859 &7774 &19288 &29633 &28718 &17546 &6811 \\
12 &0 &0 &1 &24 &359 &3178 &16578 &51799 &100810 &126013 &102743 &54469 \\

\end{tabular}
\caption{
Undirected loopless multigraphs with 2 components and a total of $E$ edges and $V$ vertices.
}
\label{tab.mk2}
\end{table}

\begin{table}
\begin{tabular}{r|rrrrrrrrrrrrrrrrrrrr}
$E\backslash V$ & 1 & 2 & 3 & 4 & 5 & 6 & 7 &8 &9 &10 &11 & 12 \\
\hline
0 &0 &0 &1\\
1 &0 &0 &0 &1\\
2 &0 &0 &0 &1 &2\\
3 &0 &0 &0 &1 &3 &4\\
4 &0 &0 &0 &1 &5 &9 &7\\
5 &0 &0 &0 &1 &6 &19 &23 &14\\
6 &0 &0 &0 &1 &9 &35 &65 &62 &29\\
7 &0 &0 &0 &1 &10 &57 &148 &214 &159 &60\\
8 &0 &0 &0 &1 &13 &89 &313 &614 &681 &421 &128\\
9 &0 &0 &0 &1 &15 &134 &601 &1577 &2374 &2148 &1104 &284\\
10 &0 &0 &0 &1 &18 &192 &1110 &3711 &7353 &8938 &6683 &3389 \\
11 &0 &0 &0 &1 &20 &269 &1938 &8225 &20752 &32692 &32639 &20712 \\
12 &0 &0 &0 &1 &24 &371 &3289 &17332 &54847 &108802 &139316 &117082 \\

\end{tabular}
\caption{
Undirected loopless multigraphs with 3 components and a total of $E$ edges and $V$ vertices.
}
\end{table}

\begin{table}
\begin{tabular}{r|rrrrrrrrrrrrrrrrrrrr}
$E\backslash V$ & 1 & 2 & 3 & 4 & 5 & 6 & 7 &8 &9 &10 &11 & 12 & 13 \\
\hline
0 &0 &0 &0 &1\\
1 &0 &0 &0 &0 &1\\
2 &0 &0 &0 &0 &1 &2\\
3 &0 &0 &0 &0 &1 &3 &4\\
4 &0 &0 &0 &0 &1 &5 &9 &8\\
5 &0 &0 &0 &0 &1 &6 &19 &24 &15\\
6 &0 &0 &0 &0 &1 &9 &35 &67 &65 &32\\
7 &0 &0 &0 &0 &1 &10 &57 &151 &221 &169 &66\\
8 &0 &0 &0 &0 &1 &13 &89 &318 &628 &711 &449 &143\\
9 &0 &0 &0 &0 &1 &15 &134 &607 &1603 &2445 &2248 &1185 &315\\
10 &0 &0 &0 &0 &1 &18 &192 &1119 &3754 &7506 &9227 &7025 &3608 \\
11 &0 &0 &0 &0 &1 &20 &269 &1949 &8294 &21049 &33426 &33790 &21799 \\
12 &0 &0 &0 &0 &1 &24 &371 &3304 &17437 &55396 &110485 &142723 &121385 \\

\end{tabular}
\caption{
Undirected loopless multigraphs with 4 components and a total of $E$ edges and $V$ vertices.
}
\label{tab.mk4}
\end{table}

\begin{table}
\begin{tabular}{r|rrrrrrrrrrrrrrrrrrrr}
$E\backslash V$ & 1 & 2 & 3 & 4 & 5 & 6 & 7 &8 &9 &10 &11 & 12 & 13 &14 \\
\hline
0 &0 &0 &0 &0 &1\\
1 &0 &0 &0 &0 &0 &1\\
2 &0 &0 &0 &0 &0 &1 &2\\
3 &0 &0 &0 &0 &0 &1 &3 &4\\
4 &0 &0 &0 &0 &0 &1 &5 &9 &8\\
5 &0 &0 &0 &0 &0 &1 &6 &19 &24 &16\\
6 &0 &0 &0 &0 &0 &1 &9 &35 &67 &66 &33\\
7 &0 &0 &0 &0 &0 &1 &10 &57 &151 &223 &172 &69\\
8 &0 &0 &0 &0 &0 &1 &13 &89 &318 &631 &718 &459 &149\\
9 &0 &0 &0 &0 &0 &1 &15 &134 &607 &1608 &2459 &2278 &1213 &330\\
10 &0 &0 &0 &0 &0 &1 &18 &192 &1119 &3761 &7533 &9299 &7126 &3690 \\
11 &0 &0 &0 &0 &0 &1 &20 &269 &1949 &8304 &21095 &33584 &34084 &22146 \\
12 &0 &0 &0 &0 &0 &1 &24 &371 &3304 &17450 &55472 &110799 &143481 &122562 \\

\end{tabular}
\caption{
Undirected loopless multigraphs with 5 components and a total of $E$ edges and $V$ vertices.
}
\label{tab.mk5}
\end{table}

\bibliographystyle{amsplain}
\bibliography{all}

\end{document}